\documentclass{article}[12pt]
\usepackage{color,times,amsmath,amsfonts,latexsym,epsfig,epsf,colordvi}
\usepackage{url,hyperref}
\usepackage[english]{babel}

\newcommand{\CC}{\mathbb {C}}
\newcommand{\RR}{\mathbb {R}}

\newcommand{\al}{\alpha }

\newcommand{\et}{\eta }
\newcommand{\la}{\lambda }

\newcommand{\bp}{\begin{pmat}}
\newcommand{\ep}{\end{pmat}}

\def\circledzero{{\begin{picture}(5.6,3.5) %
   \put(2.2,1.3){\circle{3.8}} %
   \put(1.3,.2){\normalsize 0} %
   \end{picture}}}   
\def\circledone{{\begin{picture}(5.6,3.5) 
   \put(2.2,1.3){\circle{3.8}} %
   \put(1.3,.2){\normalsize 1} %
   \end{picture}}}   
 \def\circledtwo{{\begin{picture}(5.6,3.5) %
   \put(2.2,1.3){\circle{3.8}} %
   \put(1.3,.2){\normalsize 2} %
   \end{picture}}}   
 \def\circledthree{{\begin{picture}(5.6,3.5) %
   \put(2.2,1.3){\circle{3.8}} %
   \put(1.3,.2){\normalsize 3} %
   \end{picture}}}      
   \def\circledfour{{\begin{picture}(5.6,3.5) %
   \put(2.2,1.3){\circle{3.8}} %
   \put(1.3,.2){\normalsize 4} %
   \end{picture}}}   
\def\circledfive{{\begin{picture}(5.6,3.5) %
   \put(2.2,1.3){\circle{3.8}} %
   \put(1.3,.2){\normalsize 5} %
   \end{picture}}}   
 \def\circledsix{{\begin{picture}(5.6,3.5) 
   \put(2.2,1.3){\circle{3.8}} 
   \put(1.3,.2){\normalsize 6} 
   \end{picture}}}

\author{  Frank Uhlig \thanks{Department of Mathematics and Statistics, Auburn 
University, Auburn, AL 36849-5310 \ (uhligfd@auburn.edu)}}

 \title{\vspace*{-15mm} {Zhang   Neural Networks :\\an Introduction to 
 Predictive Computations for \\ Discretized 
  time--varying Matrix Problems}}
  
\begin{document}
\date{~}
\thispagestyle{empty}
\maketitle

\vspace*{-12mm}
\hspace*{20mm}{\Large\em Dedicated -- in memoriam -- to Richard Varga \ 1928 -- 2022}
\\[14mm] 
\thispagestyle{empty}

\vspace*{-16mm}
\begin{center} {\large \bf Abstract  } \\[3mm]
\begin{minipage}{150mm}
This paper wants to increase our understanding  and computational know-how for time--varying matrix problems and Zhang Neural Networks (ZNNs). These neural networks were invented for time or single parameter--varying matrix problems around 2001 in China and almost all of their  advances have been made in and most still come from its birthplace.  Zhang Neural Network methods have become a backbone for solving discretized  sensor driven time--varying matrix problems in real-time, in theory and in on--chip applications for robots, in control theory and other engineering applications in China. They have become the method of choice for  many time--varying matrix problems that  benefit from or require  efficient, accurate and predictive real--time computations. A typical discretized Zhang Neural Network algorithm needs seven distinct steps in its initial set-up. The construction of discretized Zhang Neural Network  algorithms starts from a model equation with its associated error equation and the stipulation that the error function decrease exponentially fast. The error function differential equation is then mated with a convergent look-ahead finite difference formula to create a distinctly new multi--step style solver that predicts the future state of the system reliably from current and earlier state and solution data. Matlab codes of discretized Zhang Neural Network  algorithms for time varying matrix problems typically consist of one linear equations solve and one  recursion of already available data per time step. This makes discretized Zhang Neural network based algorithms highly competitive with ordinary differential equation initial value analytic continuation methods for function given data that are designed to work adaptively.  Discretized Zhang Neural Network  methods have  different characteristics and applicabilities than multi-step ordinary differential equations (ODEs) initial value solvers. These new time--varying matrix methods  can solve matrix--given problems from sensor data with constant sampling gaps  or from functional equations. To illustrate the adaptability of discretized Zhang Neural Networks and further the understanding of this method, this paper details the  seven step set-up process for Zhang Neural Networks and twelve separate time--varying matrix models. It supplies new codes for seven of these. Open problems are mentioned as well as detailed references to recent work on discretized Zhang Neural Networks and  time--varying matrix computations. Comparisons are given to standard non-predictive multi-step methods that use initial value problems (IVP) ODE solvers and analytic continuation methods. 
\end{minipage}\\[-1mm]
\end{center}  
\thispagestyle{empty}


\noindent{\bf AMS :}  65-02, 65-04, 65F99, 65F30, 15-04, 15A99, 15B99 \\[-7mm]

\pagestyle{myheadings}
\thispagestyle{plain}
\markboth{Frank Uhlig} {Discretized Zhang Neural Networks and time--varying Matrix Flows }

\section{Introduction }\vspace*{-1mm} 

We study and analyze  a relatively new computational approach, abbreviated occasionally as ZNN for {\em Zhang Neural Networks},  for time--varying matrix problems. The given problem's input data may come from function given matrix and vector flows $A(t)$ and $a(t)$ or from time-clocked sensor data.  Zhang Neural Networks use some of our classical notions such as derivatives, Taylor expansions (extended to not necessarily  differentiable sensor data inputs), multi-step recurrence formulas and elementary linear algebra in seven set-up steps to compute future solution data from earlier system data and  earlier solutions with ever increasing  accuracy as time progresses. \\[1mm]
In particular,  we only study discretized Zhang Neural Networks here.  
We cannot and will not  attempt to study all known variations  of Zeroing Neural Networks  (often also abbreviated by ZNN) due to the overwhelming wealth of  applications and specializations that have evolved over the last decades  with more than (estimated) 400  papers and a handful of books.  Zhang's  ZNN method has an affinity to analytic continuation methods that reformulate an algebraic system as  an ODE  initial value problem in a differential algebraic equation (DAE) and  standardly solve it by following the  solution's path via an IVP ODE initial value solver. But ZNN differs in several fundamental aspects that will be made clear in this survey. For example, Loisel and Maxwell \cite{LM18}  computed the field of values (FOV) boundary curve in 2018   via numerical continuation quickly and to high accuracy using a formulaic  expression for the FOV boundary curve of a matrix $A$ and a  single-parameter hermitean matrix flow $F(t)$ for $t \in [0, 2\pi]$. The Zhang Neural Network  method, applied to the matrix FOV problem in \cite{FUaccFOVLama20} in 2020 used the seven step  ZNN set-up and ZNN bested the FOV results of \cite{LM18}  significantly, both in accuracy and speed-wise. Unfortunately, the set-up and the workings of Zhang Neural Networks and their success has never been explained theoretically. Since the 1990s our  understandings of  analytic continuation ODE methods has been broadened and enhanced by the works of \cite{AG90, ACR94, AP98,DNM11} and others, but Zhang Neural Networks  appear to be   only adjacent to and not quite understood as part of our analytic continuation canon of adaptive multi--step DAE formulas. One obvious difference is the exponential decay of the error function $E(t)$ stipulated by $\dot E(t) = - \eta E(t)$ for $ \eta > 0$  in Zhang Neural Networks at their start versus just requiring $\dot E(t) = 0$ in analytic continuation methods. Indeed Zhang Neural Networks never solve or  try to solve their  error equation at all.\\[2mm]
Discretized ZNN methods represent a special class of Recurrent Neural Networks (RNN) that originated some 40 years ago and are intended to solve dynamical systems.  A new zeroing neural network was   proposed  by Yunong Zhang and Jun Wang in 2002, see \cite{ZW01}. 
As a graduate student at Chinese University in Hong Kong, Yunong Zhang was inspired by Gradient Neural Networks such as the Hopfield Neural Network \cite{H82} from 1982 that mimics a vector of interconnected neural nodes under time--varying neuronal inputs and has been useful in medical models and in applications to  graph theory and elsewhere. He wanted to extend Hopfield's
idea from time--varying vector algebra to more general time--varying matrix problems and his Zeroing Neural Networks, called Zhang Neural Networks or ZNN by now, were conceived in 2001 for solving dynamic parameter-varying matrix and vector problems alike.\\
Both Yunong Zhang and his Ph.D. advisor Jun Wang were unaware of ZNN's adjacency to analytic continuation ODE methods. For time--varying matrix and vector problems, Zhang and Wang's  approach starts from a global error function and an  error differential equation to achieve exponential error decay in the computed solution.\\[1mm]
  Since then,  Zhang Neural Networks  have become one mainstay for predictive time--varying matrix flow computations in the eastern engineering world. Discretized ZNN methods nowadays help with optimizing and controlling robot behavior, with autonomous vehicles, chemical plant control, image restoration, environmental sciences et cetera. They are extremely swift, accurate and robust to noise in their predictive numerical  matrix flow computations. 
  \\[1mm]
{\bf A caveat : } The term \emph{Neural Network} has had many uses.\\[1mm] 
\hspace*{6mm}\begin{minipage}{153mm}
Its origin lies in biology and medicine of the 1840s. There it refers  to the neuronal network of the brain and to the synapses of  the nervous system. In applied mathematics today  {\em neural networks} are generally associated with models  and problems that mimic or follow a  brain--like function or use nervous system like algorithms that pass information along. \\
In the computational sciences of today, assignations that use terms such as   {\em neural network} most often refer to numerical algorithms that search for relationships in parameter-dependent  data or that deal with time--varying problems. The earliest numerical use of neural networks stems from the late 1800s. Today's ever evolving numerical `neural network' methods may involve deep learning or large data and data mining. They occur in artificial neural networks, with RNNs, with continuation methods for differential equations, in homotopy methods  and in the numerical analysis of dynamical systems, as well as in  artificial intelligence (AI), in machine learning, and in image recognition and restoration,  and so forth. In each of these realizations of  'neural network' like ideas, different type algorithms are generally used for differing problems.
\end{minipage}\newpage

\hspace*{6mm}\begin{minipage}{153mm} 
{ \em Zhang Neural Networks} (ZNN) are a special zeroing neural network that differs more or less from the above. They are specifically designed to solve time--varying matrix problems and they are well suited to deal with constant sampling gap clocked sensor inputs for engineering and design applications. Unfortunately, neither
time--varying matrix problems and continuous or discretized ZNN methods are listed   in the two most recent  Mathematics Subject Classifications lists of 2010 or 2020, nor are they mentioned in Wikipedia.
\end{minipage}\\[2mm]   
In practical numerical terms,  Zhang Neural Networks and time--varying matrix flow problems  are governed by  different mathematical principles  and are subject to different quandaries than those of static matrix analysis where Wilkinson's  backward stability and error analysis are common tools and where beautiful modern static matrix principles reign.  ZNN methods can solve almost no {\em static}, i.e., {\em fixed entry} matrix problems, with the static matrix symmetrizer problem being the only known exception, see \cite{FUAZNN23}. Throughout this paper,  the term  {\em matrix flow} will describe matrices whose {\em entries} are {\em functions of time} $t$, such as the 2 by 2 dimensional matrix flow $A(t)  = \bp \sin(t^2) -t& - 3^{t-1}\\ 17.56 ~ t^{0.5}& 1/(1+t^{3.14} )\ep$.\\[0.5mm]
  Discretized ZNN processes are predictive by design. Therefore  they  require  look-ahead convergent finite difference schemes that have only rarely occurred anywhere. In stark contrast, convergent finite difference schemes but not look-ahead ones are used  in the corrector phase of multi-step  ODE solvers. \\[1mm]
Time--varying matrix computational analysis and its achievements in ZNN feel like a new and still mainly uncharted territory for Numerical Linear Algebra that is well worth studying, coding  and learning about. To begin to shed light on the foundational principles of time--varying matrix analysis is the aim of this paper.\\[2mm]
The rest of the paper is divided into two parts: Section 2 will explain the seven step set--up process of discretized ZNN in detail, see e.g., Remark 1 for some differences between Zhang Neural Networks  and analytic continuation methods for ODEs.  Section 3  lists and exemplifies a  number of  models and applied problems for time--varying matrix flows that engineers are solving, or beginning to solve, via discretized ZNN matrix algorithms. Throughout the paper we will indicate special phenomena and qualities of Zhang Neural Network  methods when appropriate.\\[-5mm]

\section{The Workings of Discretized Zhang Neural Network Methods for Para\-meter-Varying Matrix Flows}

For simplicity and by the limits of  space, given two decades of ZNN based engineering research and use of ZNN, we restrict our attention to time--varying matrix problems in discretized form throughout this paper. ZNN methods work equally well with continuous matrix  inputs by using continuous ZNN versions. Continuous and discretized ZNN  methods  have often been tested for  consistency, convergence, and stability in the Chinese literature and for their  behavior  in the presence of noise, see  the References here and  the vast literature that is available in Google.\\[1mm]
For discretized time--varying data and any matrix problem therewith, all discretized Zeroing Neural Network methods proceed using the following identical seven constructive steps -- after appropriate  start--up values have been assigned to start their iterations.
\\[1mm]
Suppose  that we are given a continuous time--varying matrix and vector model with 
time--varying functions $F$ and $G$\\[-4mm]
\begin{equation*} \hspace*{29mm}F(A(t),B(t),{x(t)}, ..) = G(t, C(t), u(t), ..) \in \RR^{m,n} \ \text{ or }  \ \CC^{m,n} \hspace*{38mm} (0)
 \end{equation*}\\[-4mm]
and a time--varying unknown vector or matrix $ x(t)$. The variables of $F$and $G$ are  compatibly sized time--varying matrices $A(t), B(t),  C(t), ...$ and time--varying vectors $u(t), ..$ that are known -- as time $t$ progresses -- at discrete equidistant time instances $t_i \in {[} t_o,t_f{]}$ for $ i\leq k$ and $k = 1, ...$ such as from sensor data. Steadily timed sensor data  is ideal for discretized ZNN.
Our task with discretized {\em predictive} Zhang Neural Networks is to find the solution $x(t_{k+1})$ of the model equation (0) accurately and in real--time from current or earlier $x(t_{..})$ values and current or earlier matrix and vector data. Note that here the 'unknown' $x(t)$ might be a concatenated vector or an  augmented  matrix $x(t)$ that may contain both, the eigenvector matrix and the associated eigenvalues for the time--varying matrix eigenvalue problem. Then the given flow matrices $A(t)$  might have  to be enlarged similarly to stay compatible with the augmented, now  'eigendata vector' $x(t)$ and likewise for other vectors or matrices $B(t), u(t)$, and so forth as needed for compatibility.\\[2mm]
Step 1 :  From a given model equation (0),  we form the {\em error function}\\[-2mm]
\begin{equation} \label{first}
E(t)_{m,n} = F(A(t),B(t),{x(t)}, ..) - G(t, C(t), u(t), ..)  \ \ \ (\stackrel{!}{=} O_{m,n} \ \text{ideally}) \end{equation}\\[-4mm]
\hspace*{12mm}which ideally should be zero, i.e., $E(t) = 0$  for all $t \in [t_o,t_f]$ if $x(t)$ solves  (0) in the desired interval.\\[1mm]
Step 2 :  Zhang Neural Networks take the derivative $\dot E(t)$ of the error function $E(t)$ and {\em stipulate} its {\em exponential \linebreak
\hspace*{12mm}decay}:\\[1mm]
\hspace*{12mm}ZNN demands that\\[-3mm]
 \begin{equation}\label{second}\dot E(t) = -\eta \ E(t)\end{equation}\\[-3mm]
\hspace*{12mm}for some fixed constant $\eta > 0$   in case of Zhang Neural Networks (ZNN).\\[1mm]
\hspace*{12mm}[ Or it demands that $$\dot E(t) = -\gamma \ {\cal F}(E(t))$$

\vspace*{-0mm}
\hspace*{12mm}for  $\gamma > 0$  and a monotonically nonlinear increasing {\em activation function} $\cal F$. Doing so changes $E(t)$\\
 \hspace*{12mm}element-wise and gives us  a  different method, called  a {\em Recurrent Neural Network} (RNN). ]\\[1mm]
The right--hand sides for ZNN and RNN methods differ subtly. Exponential error decay and thus convergence to the exact solution $x(t)$ of (\ref{second})  is automatic for both variants. Depending on the problem, different  activation functions $\cal F$ are used in the RNN  version such as  linear, power sigmoid, or hyperbolic sine functions. These can result in different and better problem suited convergence properties with RNN for highly periodic systems, see  \cite{ZYTC11}, \cite{GZ14}, or \cite{XLNSG19a} for examples.\\[0.5mm]
The exponential error decay stipulation (\ref{second}) for Zhang Neural Networks is more stringent than a PID controller would be with its stipulated error norm convergence. ZNN brings the entry--wise solution errors uniformly down to local truncation error levels,  with the exponential  decay speed depending on the value of $\eta > 0$  and it does so from any set of starting values.\\ [1mm]
As said before, in this paper we  limit our attention to discretized Zhang Neural Networks (ZNN) exclusively for simplicity.\\[1mm]
Step 3 :  Solve the exponentially decaying error equation's differential equation (\ref{second}) at time $t_k$ of Step 2\\ 
\hspace*{12mm}\underline{\bf algebraically}
for $\dot x(t_k) $ if possible. If impossible, reconsider the problem,  revise the model, and try again.\\[0.2mm]
\hspace*{12mm}(Such behavior will be encountered and mitigated in this section  and also near the end of  Section 3.) \hfill (3)\\[1mm]
Continuous or discretized ZNN never tries to solve the associated  differential error equation (\ref{second}). Throughout the ZNN process we never compute or derive  the actual computed error. The actual  errors can only be assessed by comparing the solution $x(t_i)$ with its desired quality, such as comparing $X(t_i) \cdot X(t_i)$ with $A(t_i)$ for the time--varying matrix square root problem whose relative errors are depicted in Figure 1 in Section 2. Note that ZNN cannot solve ODEs at all. ZNN is not designed for ODEs, but rather  solves time--varying matrix and vector equations.\\[0.3mm]
 The two general set-up steps 4 and 5 of ZNN below show how to  eliminate the error derivatives from the Zhang Neural Network computational process right at the start. Zhang Neural Networks do not solve or even care about the error ODE (2) at all.\\[1mm]
Step 4 $\!$: $\!\!$Select a look-ahead convergent finite difference formula for the desired local truncation error order\linebreak
 \hspace*{12mm}$O(\tau^{j+2})$ that expresses $\dot x(t_k)$ in terms of $x(t_{k+1}), x(t_k), ..., x(t_{k-(j+s)+2})$  for $ j+s $ known data points\linebreak
 \hspace*{12mm}from the table of known convergent look-ahead finite difference formulas of type {\tt j\_s} in \cite{FUfindiff19} and \cite{FUlistfindiff19}.\linebreak
 \hspace*{12mm}Here $ \tau = t_{k+1} - t_k = const$ for all $k$ is the {\em sampling gap} of the chosen discretization.\hfill (4)\\[1mm]
Step 5 :  Equate the $\dot x(t_k)$ derivative terms of  Steps 3 and 4 and thereby dispose of $\dot x(t_k)$ or the ODE problem (2)\\
\hspace*{12mm}altogether from ZNN.  \hfill (5)\\[1mm] 
Step 6 :  Solve the solution-derivative free linear equation obtained in Step 5 for $x(t_{k+1})$ and iterate in the next \linebreak \hspace*{12mm} step.  \hfill (6)\\[1mm]
Step 7 : Increase $k+1$ to $k+2$ and update all data of Step 6; then solve the updated recursion for $x(t_{k+2})$.\\ 
\hspace*{12mm}And repeat until $t_k$ reaches  the desired final time $t_f$.  \hfill (7)\\[-2mm]

For any given time--varying matrix problem the numerical Zhang Neural Network process is established in step 6. It consists of one linear equations solve (from Step 3) and one difference formula evaluation (from Step 5) per iteration step. This  predictive iteration structure of Zhang Neural Networks  with their ever decreasing errors differs  from any known analytic continuation method.\\[1mm] 
Most recently, steps 3 -- 5 above have been streamlined  in \cite{BGLZ21} by equating the Adams--Bashforth difference formula, see e.g. \cite[p. 458-460]{EMU96},  applied to the derivatives $\dot x(t_k)$ with a  convergent look--ahead difference formula from \cite{FUfindiff19} for $\dot x(t_{k})$. To achieve overall convergence, this new ZNN process requires a delicate balance between the pair of finite difference formulas and their respective decay constants  $\lambda$, and $\eta$; for further details see \cite{BGLZ21} and also \cite{DHLT93,NT11,TBD18} which explore the feasibility of $\eta$ and the sampling gap $\tau$  pairs for stability and stiffness problems with ODEs. There the feasible parameter regions are generally much smaller than what  ZNN time--varying matrix methods allow. Besides, a new  'adapted' AZNN method separates and adjusts the decay constants for different parts of ZNN individually  and thereby allows even wider $\eta$ ranges  than before and better and quicker convergence overall as well, see \cite{FUAZNN23}.\\[-2mm]

Discretized Zhang Neural Networks for  matrix problems are  highly accurate and converge quickly due to the intrinsically stipulated exponential error decay of  Step 2. These methods are still evolving and beckon new numerical analysis scrutiny and explanations which none of the 400 + applied engineering papers or any of the the handful of Zhang Neural Network books address.\\[2mm]
 The errors of Zhang Neural Nework methods have two sources: for one,  the chosen  finite difference formula's local truncation error order in Step 4 depends on  the  constant sampling gap $\tau = t_{k+1} - t_k$, and  secondly on the conditioning and rounding errors of the  linear equation solves of Step 6.  Discretized Zhang Neural Network matrix   methods are designed to give us the future solution value  $x(t_{k+1})$ accurately -- within the natural error bounds for floating--point arithmetic; and they do so  for $t_{k+1}$ immediately  after time $t_k$. At each computational step they rely only on short current and earlier equidistant sensor and solution data sequences depending on $j + s +1$ previous solution data when using a finite difference formula of type $j\_s$. \\[0.5mm]
 See formula (Pqz$i$) for example for the extended matrix eigenvalue problem several pages further down.{\\[1mm]
 Convergent  finite difference schemes have only been used in the Adams--Moulton, Adams--Bashforth, Gear, and Fehlberg  multi-step predictor--corrector formulas, see \cite[Section 17.4]{EMU96}, e.g.. We know of no other occurrence of finite difference schemes in the analytic continuation literature prior to ZNN. Discrete ZNN methods  
 can be easily transferred to on-board chip designs for driving and controlling robots and other problems once the necessary starting values for ZNN iterations have been set.  See \cite{ZGbook15} for 13 separate time--varying matrix/vector tasks, their Simulink models and circuit diagrams, as well as two chapters on fixed-base and mobile robot applications. Each chapter in \cite{ZGbook15} is referenced with 10 to 30 plus citations from the engineering literature.\\
While Zeroing Neural Networks  have been used extensively in engineering and design for two decades,  a numerical analysis of ZNN has hardly been started. Time--varying matrix numerical analysis seems to require  very different concepts than static matrix numerical analysis. Time varying matrix methods  seem to work and run according to different principles than Wilkinson's now classic backward stability and error analysis based static matrix computations. This will be made clear and clearer throughout  this introductory ZNN survey paper. \\[2mm] 
We continue our ZNN explanatory work by exemplifying the time--varying matrix eigenvalue problem $A(t)x(t) = \la(t)x(t)$ for hermitean or diagonalizable matrix flows $A(t)  \in \CC^{n,n}$ that  leads us  through the seven steps of discretized ZNN, see  \cite{BEK11} also for state of the analytic continuation based ODE  methods for the parametric matrix eigenvalue problem and \cite[Chs. 9, 10]{AG90} and \cite[Ch. 9]{AP98} for new theoretical approaches and classifications of  ODE solvers that might possibly connect discretized Zhang  Neural Networks to analytic  continuation ODE methods.   Their possible connection is not  understood  at this time and has not been researched with numerical analysis tools.\\[1mm]
 The matrix eigen--analysis  is followed by a detailed look at convergent finite difference schemes that once originated in  one-step and multi-step ODE solvers and that are  now  used predictively in discretized Zhang Neural Network methods for time--varying sensor based matrix problems.
\\[1mm] 
If $A_{n,n}$ is a diagonalizable fixed entry  matrix,  the best way to solve the static matrix eigenvalue problem $Ax = \la x$ for $A$ is to use  Francis' multi--shift implicit QR algorithm if $n \leq 11,000$ or use Krylov methods for larger sized $A$. Eigenvalues are continuous functions of the entries of $A$. Thus taking the computed eigenvalues of one flow matrix $A(t_k)$ as an approximation for the eigenvalues of $A(t_{k+1})$ might seem to suffice if the sampling gap $\tau = t_{k+1} - t_k$ is relatively small. But in practice the eigenvalues of $A(t_k)$ often  share  only a few correct leading digits with the eigenvalues of $A(t_{k+1})$. In fact the difference between any pair of respective eigenvalues of $A(t_k)$ and $A(t_{k+1})$ are generally of  size $O(\tau)$, see \cite{JZ16, YZH21, ZYGLZ21}. Hence there is need for different methods that deal more accurately  with the eigenvalues of time--varying matrix flows.\\[2mm]
By definition, for a given hermitean matrix flow $A(t)$ with $A(t) = A(t)^* \in \CC_{n,n}$  [ or for any diagonalizable matrix flow $A(t)$) ] the eigenvalue problem requires us  to compute a nonsingular matrix flow $V(t) \in \CC_{n,n}$ and a diagonal time--varying matrix flow $D(t) \in \CC_{n,n}$ so that 
$$\hspace*{55mm} A(t)V (t) = V (t)D(t) \ \text{ for all } \ t \in [t_0,t_f]. \hspace*{39mm} (0^*)$$
This is our first model equation for the time--varying matrix eigenvalue problem.\\[1mm]
 Here are  the steps for  time--varying matrix eigen--analyses when using  ZNN.\\[1mm]
Step 1 :  Create the error function  
\begin{equation*} \hspace*{50mm}E(t) = A(t)V(t) - V(t) D(t) \ \ \  (= O_{n,n} \text{ ideally.)} \hspace*{32.5mm} (1^*)
\end{equation*}
Step 2 : Stipulate exponential decay of $ E(t)$ as a function of time, i.e.,
\begin{equation*} \hspace*{60mm}\dot E(t) = -\eta \ E(t)   \hspace*{69mm} (2^*) \end{equation*}
\hspace*{12mm}for a  decay constant $\eta > 0$. \\[1mm]
\hspace*{12mm}Note : Equation $(2^*)$, written out explicitly is
\begin{eqnarray*}
\hspace*{40mm}\dot E(t) &=& \dot A(t) V(t) + A(t)\dot V(t) - \dot V(t) D(t) - V(t) \dot D(t) \nonumber \\[1mm]
& \stackrel{(*)}{=} & - \eta A(t)V(t) + \eta V(t)D(t) \ =\  - \eta E(t)\ .  \nonumber \hspace*{38mm}  
\end{eqnarray*}
\hspace*{12mm}Rearranged  with all derivatives of the unknowns $V(t)$ and $D(t)$ gathered on the left-hand side of $(*)$ :
\begin{equation*}
 \hspace*{29mm} A(t) \dot V(t) - \dot V(t) D(t) - V(t)  \dot D(t) =  - \eta A(t)V(t) + \eta V(t)D(t)  - \dot A(t) V(t)\ . \hspace*{17.6mm} (\#)
  \end{equation*}
 Unfortunately we do not know how to solve the full system eigen--equation  (\#) 
 algebraically  for the eigen--data derivative matrices $\dot V(t)$ and $\dot D(t)$ by simple matrix algebra as Step 3 asks us to do.\\[0.5mm]
  This is due to the non--commutativity of matrix products and because the unknown derivative $\dot V(t)$ appears both as a left and a right matrix factor in  the eigen--data DE  (\#). A solution that relies on Kronecker products for symmetric matrix flows $A(t) = A(t)^T$ is available in \cite{ZHYLLQ20} and we will follow the Kronecker product  route later   when dealing with square roots of time--varying matrix flows in subparts ({\bf VII}) and ({\bf VII start-up}) in this section, as well as  when solving time--varying classic matrix equations via ZNN in subpart ({\bf IX}) of Section 3.\\[-1mm]
  \hspace*{70mm} \underline{\hspace*{10mm}}\\[1.5mm]
 Now we revise our matrix eigen--data model and  restart the whole process anew. To overcome the above dilemma, we separate the global time--varying matrix eigenvalue problem for $A_{n,n}(t)$ into $n$ eigenvalue problems 
 $$\hspace*{52mm}A(t)x_i(t) = \la_i(t)x_i(t) \ \text{ with } \ i = 1,...,n \hspace*{44.5mm} (0i)$$ 
 that can be solved for one eigenvector $x_i(t)$ and one eigenvalue $\la_i(t)$ at a time as follows.\\[1mm]
 Step 1 : The error function $(0i)$ in vector form is
 \begin{equation*} \hspace*{47mm}e(t) = A(t)x_i(t) - \la_i(t) x_i(t) \ \ \  (= o_{n} \in \CC^n \ \ \text{ ideally)} \hspace*{34mm} (1i)
\end{equation*}
Step 2 : We demand exponential decay of $ e(t)$ as a function of time, i.e., 
\begin{equation*} \hspace*{58mm}\dot e(t) = -\eta~ e(t) \hspace*{75.6mm} (2i) \end{equation*}
\hspace*{12mm}for a  decay constant $\eta > 0$. \\[1mm]
\hspace*{12mm}Equation $(2i)$, written out explicitly, now becomes
\begin{eqnarray*}
\hspace*{36mm}\dot e(t) &=& \dot A(t) x_i(t) + A(t)\dot x_i(t) - \dot \la_i(t)  x_i(t)  -  \la_i(t) \dot x_i(t)\\[1mm]
& \stackrel{(*)}{=} & - \eta~ A(t)x_i(t) + \eta~ \la_i (t) x_i(t)\ =\  - \eta~ e(t).   \hspace*{42mm} 
\end{eqnarray*}
\hspace*{12mm}Rearranged,  with the derivatives of $x_i(t)$ and $\la_i(t)$ gathered on the left--hand side of $(*)$ :
\begin{equation*}
\hspace*{21mm}  A(t) \dot x_i(t) - \dot \la_i(t)  x_i(t) -   \la_i(t)  \dot x_i(t)=  - \eta~ A(t)x_i(t) + \eta~ \la_i(t)  x_i(t)  - \dot A(t) x_i(t)\ . \hspace*{23.2mm} 
 \end{equation*} 
\hspace*{12mm}Combining the $\dot x_i$ derivative terms gives us
$$ \hspace*{24mm} ( A(t) - \la_i(t) I_n) \dot x_i(t) - \dot \la_i(t)  x_i(t)  =    (- \eta~ (A(t) - \la_i(t) I_n)   - \dot A(t)) ~ x_i(t)\ . \hspace*{28.4mm} 
 $$
 For each $i = 1,...,n$  this equation  is a  differential equation in the  unknown eigenvector ${x_i(t)} \in \CC^n$ and the  unknown eigenvalue ${\la_i(t)} \in \CC$.  We concatenate $x_i(t)$ and $\la_i(t)$ in $$z_i(t) = \bp x_i(t)\\\la_i(t) \ep \in \CC^{n+1}$$ and obtain the following matrix/vector DE for the  eigenvector $x_i(t)$ and its associated eigenvalue $\la_i (t)$ and each $ i = 1, ..., n$, namely
 $$\hspace*{14.5mm} \begin{pmat} A(t) - \la_i(t) I_n & -x_i(t) \end{pmat}_{n,n+1}
 \begin{pmat} \dot x_i(t)\\ \dot \la_i(t) \end{pmat}
    =   \left( - \eta~ (A(t) - \la_i(t) I_n)   - \dot A(t) \right) x_i(t) \in \CC^n   \  \hspace*{4mm} \text{(Az}i)
$$
where the augmented system matrix on the left--hand side of formula (Az$i$) has dimensions $n$ by $n+1$ if $A$ is $n$ by $n$.\\[1mm]
Since each matrix eigenvector  defines an  invariant 1--dimensional subspace we must ensure that   the computed eigenvectors $x_i(t)$ of $A(t)$ do not grow infinitely small or infinitely large in their ZNN  iterations.  Thus we require  that the computed  eigenvectors   attain unit  length asymptotically by introducing the additional error function $e_2(t) = x^*_i(t) x_i(t) - 1$. Stipulating exponential decay for $e_2$ leads to 
$$ \dot e_2(t) = 2 x^*_i(t) \dot x_i(t) = -\mu~ ( x^*_i(t) x_i(t) - 1) = -\mu~ e_2(t)\  $$\\[-5mm]
or 
\vspace*{-3mm}
$$\hspace*{50.5mm} -x^*_i(t) \dot x_i(t) = \mu/2~ ( x^*_i(t) x_i(t) - 1) \hspace*{51mm} \text{(e}_2i)$$ 
for a second decay constant $\mu > 0$.
If we set $\mu = 2 \eta$, place equation (e$_2i$) below the last row of the $n$ by $n+1$ system matrix of equation (Az$i$), and extend its right--hand side vector by the right--hand side entry in (e$_2i$), we obtain an $n+1$ by $n+1$ time--varying system of DEs (with a hermitean system matrix if $A(t)$ is hermitean). I.e., \\
$$\hspace*{16mm} \begin{pmat}
 A(t) - \la_i(t) I_n & -x_i(t)\\
 -x_i^*(t) & 0 
 \end{pmat}
 \begin{pmat} \dot x_i(t)\\ \dot \la_i(t) \end{pmat}
    =   \begin{pmat} ( - \eta~ (A(t) - \la_i(t) I_n)   - \dot A(t)) x_i(t) \\ \eta~ (x_i^*(t)x_i(t) -1) \end{pmat} \ .\hspace*{6.9mm} \text{(Pqz}i)   
$$
\hspace*{12mm}Next we set
\begin{equation*}
\begin{array}{c}
 P(t_k) = \begin{pmat}
 A(t_k) - \la_i(t_k) I_n & -x_i(t_k)\\
 - x_i^*(t_k) & 0 
 \end{pmat} \in \CC_{n+1,n+1} , \  \ \ z(t_k)= \begin{pmat}  x_i(t_k)\\  \la_i(t_k) \end{pmat} \in \CC^{n+1} \ , \\[7mm]
 \text{ and }  \ q(t_k) = \begin{pmat} ( - \eta~ (A(t_k) - \la_i(t_k) I_n)   - \dot A(t_k)) x_i(t_k) \\ \eta~ (x_i^*(t_k)x_i(t_k) -1)  \end{pmat} \in \CC^{n+1} 
 \end{array}
 \end{equation*}
 \hspace*{12mm}for all discretized  times $t = t_k$.  And we have completed Step 3 of the ZNN set-up process.\\[2mm] 
 Step 3 : Our model  $(0i)$ for the $i$th eigenvalue equation of $A(t_k)$ has been transformed into  the mass matrix/vector\\   
 \hspace*{12mm}differential equation  \\[-3mm]
 \begin{equation*}
 \hspace*{32mm}P(t_k) \dot z(t_k) = q(t_k) \ \ \ \text{ or } \ \ \  \dot z(t_k) = P(t_k) \backslash q(t_k) \ \ \text{ in Matlab notation}.\hspace*{25mm} (3i)   
 \end{equation*}
 Note that at time $t_k$ the mass matrix $P(t_k) $ contains the current input data $A(t_k)$ and currently computed eigen--data $\la_i (t_k)$ and $x_i(t_k)$ combined in the vector  $z(t_k)$, while the right-hand side vector $q(t_k)$ contains  $A(t_k)$ as well as its first derivative $\dot A(t_k)$, the decay constant $\eta$, and the computed eigen--data at time $t_k$.\\[1mm] 
 The system matrix $P(t_k)$ and the right--hand side vector $q(t_k)$  in formula (3i) differ greatly from the simpler differentiation formula for the straight DAE used in \cite{LM18} where $\dot E(t) =0$ or $\eta =0 $ was assumed.\\[1mm]
 Step 4 : Now we  choose the  following convergent look-ahead finite 5-IFD (five Instance Finite Difference)  formula\\[.5mm]
\hspace*{12mm}of type {\tt j\_s} = {\tt 2\_3} with global truncation error order $O(\tau^3)$ from the list in \cite{FUlistfindiff19}  for $\dot z_k$ : \\[-1mm]
 \begin{equation*}
\hspace*{40mm}\dot z_k = \dfrac{8z_{k+1}+z_k -6z_{k-1} - 5 z_{k-2} + 2 z_{k-3}}{18 \tau} \in \CC^{n+1} \ .\hspace*{36mm}(4i) 
\end{equation*}

\vspace*{0mm}
Step 5 : Equating the different expressions for $18 \tau \dot z_k$ in $(4i)$ and $(3i)$ (from steps 4 and 3 above) we have\\[-3mm]
$$\hspace*{30mm} 18\tau \cdot \dot z_k =8z_{k+1}+z_k -6z_{k-1} - 5 z_{k-2} + 2 z_{k-3} \stackrel{(*)}{=} 18 \tau \cdot (P \backslash q) = 18\tau \cdot \dot z_k\hspace*{20mm} (5i)
$$

\vspace*{-1mm}
\hspace*{12mm}with local truncation error order $O(\tau^4)$  due to the multiplication of both $(3i)$ and $(4i)$ by $18 \tau$. \\[2mm]
Step 6 : Here we solve the inner equation $(*)$ in $(5i)$  for $z_{k+1}$ and obtain the discretized look--ahead ZNN iteration\\[.5mm]
\hspace*{12mm}formula 
$$ \hspace*{32mm} z_{k+1} = \dfrac{9}{4} \tau (P(t_k) \backslash q(t_k)) - \dfrac{1}{8} z_k + \dfrac{3}{4} z_{k-1} + \dfrac{5}{8} z_{k-2} - \dfrac{1}{4} z_{k-3} \ \in \CC^{n+1} \hspace*{25mm} (6i)
$$
\hspace*{12mm}that is comprised of  a linear equations part and a recursion part and has local truncation error order $O(\tau^4)$.\\[2mm]
Step 7 : Iterate to predict the eigendata vector  $z_{k+2}$  for $A(t_{k+2})$ from earlier eigen and system  data at times $t_{\tilde j}$ with\\
\hspace*{12mm}$\tilde j \leq k+1$ and repeat. \hfill (7i)\\[2mm]
The final formula $(6i)$ of ZNN contains the computational formula for future eigen--data with near unit eigenvectors in the top $n$ entries of $z(t_{k+1})$ and the eigenvalue appended below.  Only a  mathematical  formula of type $(6i)$ needs to be derived and implemented in code for any other matrix model problem. In our specific case and for any other time--varying matrix problem  all  entries in $(6i)$ have been computed earlier as eigen--data for times $t_{\tilde j}$ with $\tilde j \leq k$, except for the system or sensor input $A(t_k)$ and $\dot A(t_k)$. The derivative  $\dot A(t_k)$ is best computed via a high error order derivative formula from previous $\dot A(t_{\tilde j})$ with $\tilde j \leq k$.\\[2mm]
The  computer  code lines that perform the actual math for ZNN iteration steps for a single eigenvalues from time $t = t_k$ to $t_{k+1}$ and $k = 1, 2,3,...$ are listed below. There $ze$ denotes one eigenvalue of $Al = A(t_k)$ and $zs$ is its associated eigenvector. $Zj$ contains the relevant set of earlier eigen--data for $A(t_j)$ with $j \leq k$  and $Adot$ is an approximation for $\dot A(t_k)$. Finally  $eta, tau$ and $taucoef\!f$ are chosen for  the desired error order finite difference formula and its characteristic polynomial, respectively.\\[-7mm]
\begin{verbatim}
              .              .             .
              .              .             .
    Al(logicIn) = diag(Al) - ze*ones(n,1); % Al = A(tk) - ze*In
    P = [Al,-zs;-zs',0];                   % P is generally not hermitean
    q = taucoeff*tau*[(eta*Al + Adot)*zs; -1.5*eta*(zs'*zs-1)]; % rh side
    X = linsolve(P,q);                     % solve a linear equation 
    Znew = -(X + Zj*polyrest);             % New eigendata at t_{k+1}
    ZN(:,jj) = Znew;                       % Extend the known eigendata
              .              .            .
              .              .            .
\end{verbatim}

\vspace*{-30.5mm} \hspace*{6mm}\rule{1mm}{15mm}\\[14mm]
The four central code lines above that are  barred along the left edge  express the computational essence of Step 6 for the matrix eigen--problem in ZNN. Only there is any math performed, the rest of the program code are input reads and output saves and preparations. After Step 6 we  store the new data and repeat these 4 code lines with peripherals for $t_{k+2}$ until we are done  with  one eigenvalue at $t = t_f$. Then we repeat the same code for the next eigenvalue of a hermitean or diagonalizable matrix flow $A(t)$.\\[2mm]
What is actually computed   in  formula $(6i)$ for time--varying matrix eigen--problems  in ZNN? To explain we recall some convergent finite difference formula theory next.\\[1mm]
It is well known that the characteristic polynomial coefficients $a_k$ of a    
  finite difference scheme must add up to zero for convergence, see \cite[Section 17]{EMU96} e.g.. Thus for a look--ahead and convergent finite difference formula \\[-2.5mm]
\begin{equation*} z(t_{k+1}) + \al_k z(t_k) + \al_{k-1} z(t_{k-1}) + \ \cdots \ + \al_{k-\ell} z(t_{k-\ell}) \end{equation*}\\[-4mm]
and its  { characteristic polynomial} \\[-2.5mm]
$$ p(x) = x^{k+1} + \al_k x^k + \al_{k-1} x^{k-1} + \ \cdots \ + \al_{k-\ell} x^{k-\ell} $$
we must have that $p(1) = 1 + \al_k + \al_{k-1} + ... + \al_{k-\ell} = 0$. Plugging $z_{..}$ into the 5-IFD difference formula $(4i)$, we realize that  in formula $(6i)$  \\[-2mm]
$$\hspace*{34mm}z(t_{k+1} ) +\dfrac{1}{8} z(t_k) - \dfrac{3}{4}  z(t_{k-1}) - \dfrac{5}{8} z(t_{k-2})  + \dfrac{1}{4} z(t_{k-3}) \approx o_{n+1} \hspace*{31.2mm} (8i)$$\\[-3mm] 
 due to unavoidable truncation and rounding errors. Thus asymptotically and with the stipulated exponential error decay\\[-2mm] 
  $$\hspace*{39mm} z(t_{k+1} ) \approx -\dfrac{1}{8} z(t_k) + \dfrac{3}{4}  z(t_{k-1}) +\dfrac{5}{8} z(t_{k-2})  - \dfrac{1}{4} z(t_{k-3}) \hspace*{36mm} (9i)$$\\[-3mm] 
  with local truncation error of order $O(\tau^4)$ for the chosen 5-IFD {\tt j\_s} = {\tt 2\_3}  difference formula in $(4i)$.\\[1.5mm]
 In step 6 of the above ZNN process, formula $(6i)$  
 $$  z_{k+1} =  \dfrac{9}{4} \tau (P(t_k) \backslash q(t_k)) - \dfrac{1}{8} z_k + \dfrac{3}{4} z_{k-1} + \dfrac{5}{8} z_{k-2} - \dfrac{1}{4} z_{k-3} \ \in \CC^{n+1} 
$$
splits the 5-IFD formula  $(8i)$ into two nearly equal parts that become ever closer to each other in $(9i)$ due to the nature of our convergent finite difference schemes. 
The first term \ $\dfrac{9}{4} \tau (P(t_k) \backslash q(t_k)) $  in $(6i)$ adjusts the predicted value of $z(t_{k+1})$ only slightly  according to the current system data inputs while the remaining finite difference formula term $$- \dfrac{1}{8} z_k + \dfrac{3}{4} z_{k-1} + \dfrac{5}{8} z_{k-2} - \dfrac{1}{4} z_{k-3}$$ in $(6i)$ has eventually a nearly identical magnitude as  the solution vector $z_{k+1}$ at  time $t_{k+1}$. \\[1mm]
And indeed for the time--varying matrix square root problem in {\bf (VII)}, our test  flow matrix $A(t)$  in Figure 1 increases in norm to around 10,000 after 6 minutes of simulation and  the norm of the solution square root matrix flow $X(t)$ hovers around 100 while the magnitude of the first linear equations solution term of $(6i)$ is around $10^{-2}$, i.e., the two terms of the ZNN iteration $(6i)$ differ in magnitude by a factor of around $10^4$, a disparity  in magnitude that we should expect from the above analysis.
This behavior of ZNN matrix methods is exemplified by the error graph for  time--varying matrix square root computations  in Figure 1 below that will be discussed further in Section 3, part {\bf (VII)}. The 'wiggles' in the error curve of Figure 1 after the initial decay phase represent the relatively small input data adjustments that are made by the linear solve term of ZNN. For more details and data on the magnitude disparity see \cite[Fig 6, p.173]{FUAZNN23}.\\[1mm]
Figure 1 was computed via simple Euler steps from a random entry matrix  start-up matrix  in\linebreak  
 {\tt tvMatrSquareRootwEulerStartv.m}, see Section 3 {\bf  (VII start-up)}. The main ZNN iterations for Figure 1 have  used a 9-IFD formula of type {\tt 4\_5} with local truncation order 6.\\[-4mm]
\begin{center}
 \hspace*{12 mm}\includegraphics[width=100mm]{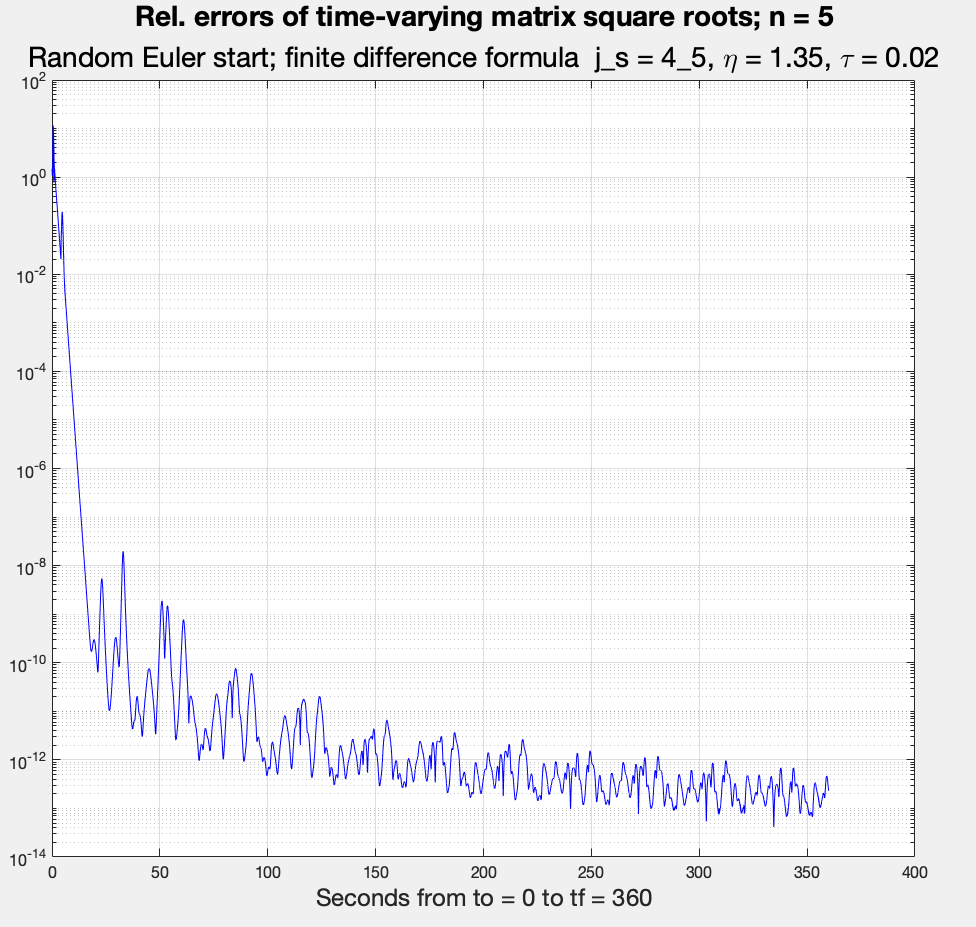}\\[-1mm]
 Figure 1 : Typical relative error output for {\tt  tvMatrSquareRootwEulerStartv.m} in Section 3 {\bf (VII)}
\end{center}
The 5-IFD formula in equation (4i) above is of type {\tt j\_s = 2\_3}  and   in discretized ZNN its  local truncation error order is relatively low at $O(\tau^4)$ as $j+2 = 2+2 = 4$. In tests we prefer to  use a 9-IFD of type  {\tt 4\_5}. To start a discretized ZNN iteration process with a look-ahead convergent finite difference formula of type {\tt j\_s} from the list in  \cite{FUlistfindiff19} requires $j+s$ known starting values. For  time--varying matrix eigenvalue problems that are given by function inputs for $A(t_k)$ we generally use Francis QR to generate the $j+s$ start--up eigen--data set, then iterate via discretized ZNN. And throughout the discrete iteration process from $t = t_o$ to $t = t_f$ we need to keep  only the most recently computed  $j+s$  points of data to compute the eigen--data predictively for  $t = t_{k+1}$. \\
MATLAB codes for several time--varying matrix  eigenvalue computations via discretized ZNN are  available at \cite{FUZNNEigV}.\\[2mm]
Next we study how to construct general look-ahead and convergent finite difference schemes of arbitrary truncation error orders $O(\tau^p)$ with $p \leq 8$ for use in discretized ZNN time--varying matrix methods. We start  from random entry  seed vectors and can use Taylor polynomials and elementary linear algebra to construct look-ahead finite difference formulas of any order that may or -- most likely -- may not be convergent. The constructive first step of finding look-ahead finite difference formulas is followed by a second, an optimization procedure to find look-ahead \underline{and} convergent finite difference formulas of the desired error order. This second non-linear part  may  not always succeed as we shall explain later.\\[1mm] 
 Consider a discrete time--varying state vector $x_k = x(t_k) = x(t_o+k \cdot \tau)$ for a constant sampling gap $\tau$ and $k = 0,1,2,...$ and write out $r = \ell + 1$ explicit Taylor expansions of degree $j$ for $x_{k+1}, x_{k-1}, ..., x_{k-\ell}$ about $x_k$.\\[1mm]
Each Taylor expansion in our scheme  will contain $j+2$ terms on the right hand side, namely $j$ derivative terms, a term for $x_k$, and one for the error $O(\tau^{j+1})$. Each right hand side's under- and over-braced $j-1$ 'column'  terms in the scheme displayed below contain  products of identical powers of $\tau$  and identical higher derivatives of $x_k$ which -- combined in column vector form, we will call $taudx$. Our aim is to find a linear combination of these $r = \ell +1$ equations for which the under- and over-braced sums vanish for all possible higher derivatives of the solution $x(t)$ and all sampling gaps $\tau$. If we are able to do so,  then we can express $x_{k+1}$ in terms of $x_l$ for $l = k,k-1,...,k-\ell$, $\dot x_k$, and $\tau$ with a local truncation error of order $O(\tau^{j+1})$ as a linear combination of the depicted $r=\ell + 1$  Taylor equations (10) through (15) below. Note that the first Taylor expansion (10) for $x(t_{k+1})$ is unusual, being  look-ahead. This has never been used in standard  ODE schemes. Equation (10) sets up the predictive behavior of ZNN.\\[2mm]
\setcounter{equation}{9}
 \hspace*{78mm} \raisebox{0.5mm}{$j-1$ terms}\\[-6.4mm]  
\begin{eqnarray}
x_{k+1} & =  & x_k ~ +~  \tau \dot{x}_k \ \overbrace { + \ \ \dfrac{\tau^2}{2!} ~ \ddot{x}_k \ \ ~ + ~ \  \ \dfrac{\tau^3}{3!} ~ \overset{\dots}{x}_k  ~ \ \  ... \ \ \ \ \ \ \ \ \ \ + \  \dfrac{\tau^{j}}{j!} ~ \overset{j}{\dot{x}}_k \ \ \ \ \ \   } \  \  \ + \ ~O(\tau^{j+1})\\[1mm]
x_{k-1} & =  & x_k ~ - ~ \tau \dot{x}_k \ + \ \dfrac{\tau^2}{2!} ~ \ddot{x}_k \ \ ~ -  ~ \ \ \dfrac{\tau^3}{3!} ~ \overset{\dots}{x}_k  ~ \ \ ...  \    \ \ \ +   (-1)^{j}\dfrac{\tau^{j}}{j!} ~ \overset{j}{\dot{x}}_{k} \ \ \ \ \ ~  + \  O(\tau^{j+1})\\[1mm]
x_{k-2} & =  & x_k - 2\tau \dot{x}_k \ + \dfrac{(2\tau)^2}{2!} ~ \ddot{x}_k -  \dfrac{(2\tau)^3}{3!} ~ \overset{\dots}{x}_k  ~ ... ~ + (-1)^{j} \dfrac{(2\tau)^{j}}{j!} ~ \overset{j}{\dot{x}}_{k} \ \ \ + \ O(\tau^{j+1})\\[1mm]
x_{k-3} & =  & x_k - 3\tau \dot{x}_k \ + \dfrac{(3\tau)^2}{2!} ~ \ddot{x}_k -  \dfrac{(3\tau)^3}{3!} ~ \overset{\dots}{x}_k ~  ...  ~ + (-1)^{j} \dfrac{(3\tau)^{j}}{j!} ~ \overset{j}{\dot{x}}_{k} \ \ \ + \ O(\tau^{j+1})\\[2mm]
& \vdots& \hspace*{53mm} \vdots \\[2mm]
x_{k-\ell} & =  & x_k - \ell\tau \dot{x}_k \ \ \underbrace { + \ \dfrac{(\ell\tau)^2}{2!} ~ \ddot{x}_k -  \dfrac{(\ell\tau)^3}{3!} ~ \overset{\dots}{x}_k ~  ...  \ + (-1)^{j} \dfrac{(\ell\tau)^{j}}{j!} ~ \overset{j}{\dot{x}}_{k}   } \ \ ~ + ~ O \ (\tau^{j+1})
\end{eqnarray}

\vspace*{-2.5mm}
\hspace*{78mm} $j-1$ terms\\[1mm]
 The 'rational number' factors in  the 'braced  $j-1$ columns' on the right--hand side of  equations (10) through (15) are collected in the $r = \ell + 1$ by $j-1$ matrix $\cal A$ :
\begin{equation}\label{ratA} {\cal A}_{r,j-1} = \begin{pmat} \dfrac{1}{2!} & \dfrac{1}{3!} & \dfrac{1}{4!} & \cdots & \cdots & \dfrac{1}{j!}\\[5mm]
   \dfrac{1}{2!} & -\dfrac{1}{3!} & \dfrac{1}{4!} & \cdots & \cdots &  (-1)^{j}~ \dfrac{1}{j!}\\[5mm]
   \dfrac{2^2}{2!} & - \dfrac{2^3}{3!} & \dfrac{2^4}{4!} & \cdots & \cdots & (-1)^{j}~ \dfrac{2^{j}}{j!}\\[3mm]
     \vdots & \vdots & \vdots && & \vdots\\[2mm]
  \dfrac{\ell^2}{2!} & - \dfrac{\ell^3}{3!} & \dfrac{\ell^4}{4!} & \cdots & \cdots & (-1)^{j}~ \dfrac{\ell^{j}}{j!} 
\end{pmat} _{r,j-1}  .
\end{equation}
Now the over-- and under--braced expressions in    equations (10) through (15)  above have the  matrix times vector product form\\[-3mm]
\begin{equation}\label{ADx}
{\cal A} \cdot  taudx = \begin{pmat} \dfrac{1}{2!} & \dfrac{1}{3!} & \dfrac{1}{4!} & \cdots & \cdots & \dfrac{1}{j!}\\[5mm]
   \dfrac{1}{2!} & -\dfrac{1}{3!} & \dfrac{1}{4!} & \cdots & \cdots &  (-1)^{j}~ \dfrac{1}{j!}\\[5mm]
   \dfrac{2^2}{2!} & - \dfrac{2^3}{3!} & \dfrac{2^4}{4!} & \cdots & \cdots & (-1)^{j}~ \dfrac{2^{j}}{j!}\\[3mm]
     \vdots & \vdots & \vdots && & \vdots\\[2mm]
  \dfrac{\ell^2}{2!} & - \dfrac{\ell^3}{3!} & \dfrac{\ell^4}{4!} & \cdots & \cdots & (-1)^{j}~ \dfrac{\ell^{j}}{j!} 
\end{pmat}_{r,j-1} \hspace*{-6.7mm} 
\cdot
\begin{pmat}
\tau^2 \ \ddot{x}_k \\[1mm] \tau^3 \ \overset{\dots}{x}_k \\[1mm] \tau^4 \ \overset{4}{\dot{x}}_k \\[1mm] \vdots \\[1mm]   \vdots \\[1mm] \tau^{j} \ \overset{j}{\dot{x}}_{k}
\end{pmat}_{j-1,1} 
\end{equation}
where the $j-1$-dimensional column vector $taudx$ contains the increasing powers of $\tau $ multiplied by the respective higher derivatives of $x_k$ that appear in the Taylor expansions (10) to (15).\\[1mm]
Note that for any nonzero  left kernel row vector $y \in \RR^{r}$ of ${\cal A}_{r,j-1}$ with $y \cdot {\cal A} = o_{1,j-1}$ we have

\vspace*{-3mm}  
$$y_{1,r} \cdot {\cal A}_{r, j-1} \cdot taudx_{j-1,1} = o_{1,j-1} \cdot taudx_{j-1,1} = 0\in \RR \ ,
$$  
no matter what the higher derivatives of $x(t)$ at $t_k$ are. Clearly we can zero  all under-    and over--braced terms in equations (14) to (19)  if ${\cal A}_{r,j-1} $  has a  nontrivial left kernel. This is certainly the case if $\cal A$ has  more rows than columns, i.e., when   \ $r=\ell+1 >j-1$. A non-zero left kernel vector $w$ of $\cal A$ can then be easily found via Matlab, see \cite{FUfindiff19} for details.
The linear combination of  the equations (10) through (15) with the coefficients of  $w \neq o$ creates a  predictive recurrence relation  for $x_{k+1}$ in terms of $x_k, x_{k-1}, ..., x_{k-\ell}$ and $\dot x_k$ with local truncation error order $O(\tau^{j+2})$ as desired. Solving this recurrence equation for $\dot x_k$ gives us formula $(4i)$ in Step 4. Then multiplying by a multiple of $\tau$ in Step 5 increases the resulting formula's local truncation error order from $O(\tau^{j+1})$ to $O(\tau^{j+2})$ in equation $(5i)$.\\[2mm]
 Thus far we have tacitly assumed that the solution $x(t)$ is sufficiently often differentiable for high order Taylor expansions. This can rarely be the case in real world applications, but -- fortunately -- discretized computational ZNN works very well with errors as predicted by Taylor even for discontinuous and limited random sensor failure inputs and consequent non-differentiable $x(t)$. The reason is a mystery. We do not  know of any non-differentiable Taylor expansion theory. What this might be is a challenging open question for time--varying numerical matrix analysis.\\[1mm] 
A recursion formula's  characteristic polynomial determines its convergence and thus its suitability for discretized ZNN methods. More specifically, the convergence of finite difference formulas and recurrence relations like ours hinges on the lay of the roots of their  associated characteristic polynomials in the complex plane. This  is  well known for multi-step  formulas and also applies to processes such as discretized ZNN recurrences.  Convergence requires that all roots of the formula's characteristic polynomial  lie inside or on the unit disk in $\CC$ with no repeated roots allowed on the unit circle, see \cite[Sect. 17.6.3]{EMU96} e.g..\\[2mm]
Finding convergent {\em and} look--ahead finite difference formulas from look-ahead ones is a non-linear problem. In \cite{FUfindiff19} we have approached this problem by  minimizing  the   maximal modulus root of  'look--ahead' characteristic polynomials to below $1 + eps$ for a very small threshold $0 \approx eps\geq 0$ so that they become numerically and practically convergent  while the polynomials' coefficient vectors $w$ lie in the left kernel of ${\cal A}_{r,j-1}$.\\[1mm]
The set of look-ahead characteristic polynomials is not a  subspace since sums of such polynomials may or -- most often -- may not represent look-ahead finite difference schemes. Hence we must  search  indirectly  in a  neighborhood of the starting seed $y\in \RR^{r-j+1}$ for   look--ahead  characteristic polynomials with proper minimal maximal magnitude roots. 
For this  indirect root minimization process we have used  Matlab's multi--dimensional minimizer function {\tt fminsearch.m}  that uses the  Nelder--Mead downhill simplex method, see \cite{NM65, LRWW98}.  It mimics the method of  steepest descent and searches for local minima via multiple function evaluations  without using derivatives until it has either found a look--ahead seed with an associated  characteristic polynomial that is convergent and for which its coefficient vector remains in the left kernel of the associated ${\cal A}_{r,j}$ matrix, or there is no such convergent formula from the chosen look-ahead seed. \\[1mm]
Our two part look--ahead and convergent difference formula  finding algorithm has computed many convergent and look--ahead finite difference schemes for discretized ZNN of all types {\tt j\_s}  with $1 \leq j \leq 6$ and $j \leq s \leq j+3$ with local truncation  error orders between $O(\tau^3)$ and $O(\tau^8)$.  Convergent look--ahead finite difference formulas  were unavailable before  for ZNN use with error orders above $O(\tau^4)$. And different low error order formulas had only been used  before in the corrector phase of predictor-corrector ODE solvers.\\[1mm]
In our experiments with trying to find convergent and look--ahead discretization formulas of type {\tt j\_s} where $s = r-j$ we have never succeeded when $1 \leq s = r-j  <j$. Success always occurred for $s = j$ and better success when $s = j+1$ or $s=j+2$. It is obvious that for $s= r-j= 1$ and any seed $y \in \RR^1$ there is only one normalized 
look--ahead discretization formula {\tt j\_1} and it appears to  never be  convergent. For convergence we seemingly need more freedom in our seed space $\RR^{r -j}$ than there is in one dimension or  even in less than $j$--dimensional space.\\[1mm]
A different method to find convergent look--ahead finite difference formulas starting from the same Taylor expansions in (10) to (15) above has been derived in \cite{YZH21} for various IFD formulas with  varying truncation error orders from 2 through 6 and rational coefficients. These two methods have not been compared or their differences studied thus far.
\newpage

{\bf Remark 1 : } \\[2mm]
\hspace*{4mm}\begin{minipage}{154mm}{Discretized Zhang Neural Networks  were  originally designed for sensor based time--varying matrix problems. They solve sensor given matrix flow  problems accurately in real--time and for on--line chip applications. They are being used  extensively in this way today. Yet discretized matrix Zhang Neural Networks   are partially built on one differential equation, the error function DE (2), and they use custom multi--step finite difference formulas as our  centuries old numerical initial value ODE solvers do in analytic continuation algorithms.\\[1mm]
 What are their differences? I.e., can discretized ZNN be used, somehow  understood, or   interpreted as part of an analytic continuation method for ODEs.  Or is the reverse possible, i.e., can numerical IVP ODE solvers be used successfully  for function based  parameter--varying matrix problems?\\[1mm]
 Only the latter seems possible, see \cite{LM18} for example. There the accuracy and speed comparison 
 of the much more difficult to evaluate  field of values boundary curve equation is bettered by an analytic continuation  method that integrates the computationally simpler derivative more accurately and quickly.\\[0.5mm] 
 Yet for time--varying, both  function based and sensor based, matrix flow problems Zhang Neural Networks algorithms are seemingly a breed of their own. \\[1mm]
  In Section 3 (V) we use them to solve time--varying linear systems of nonlinear equations via Lagrange matrix multipliers and  linearly constrained nonlinear optimization.\\[1mm]
  Here, however,  we deliberately focus on the seven step design of  Zhang Neural Networks for discretized time--varying  matrix flow problems.\\[1mm]
 Given an initial value function based ODE problem\\[-3mm]
$$ y'(t) = f(t,y(t)) \ \  \text { with } \ \ y(t_o) = y_o \ \text{ and } t \in [a,b] ,$$
and if we formally integrate
$$ \int_{t_i}^{t_{i+1}} y'(t) \ dt = \int_{t_i}^{t_{i+1}}  f(t,y(t)) \ dt  \ \text{ with } t_i \geq a \text{ and } t_{i+1} \leq b 
$$
we obtain the one-step look-ahead Euler-type formula\\[-2mm] 
\vspace*{-0mm} 
$$y(t_{i+1}) = y(t_i) + \int_{t_i}^{t_{i+1}}  f(t,y(t)) \ dt $$\\[-2mm]
for all feasible $i$. The aim of IVP ODE solvers is to compute an accurate  approximation for each partial integral of $f$ from $t_i$ to $t_{i+1}$ and then approximate the  value of $y$ at $t_{i+1}$ by adding $\int_{t_i}^{t_{i+1}}  f(t,y(t)) dt $ to its previously computed value $y(t_i)$.
One-step methods use ever more accurate integration formulas to achieve an accurate table of values for $y$,  provided that  an antiderivative of $y'$ is not known, cannot be found in our antiderivative tables, and cannot be computed via symbolic integration software.\\[0.5mm]
 The multi-step formulas of Adams--Bashforth,  Adams--Moulton and others  use several earlier computed $y$
values and integration formulas in tandem in a two-process method of prediction followed by a correction step. The Fehlberg multi--step corrector formula \cite[p.464]{EMU96}, \cite{F61} uses three previous $y$ values and five earlier right--hand side function $f$ evaluations while the Gear corrector formulas \cite[p. 480]{EMU96}, \cite{G71/1}
  use up to six previously  computed  $y$ values and just one function evaluation to obtain $y(t_{i+1})$ from data at or before time $t_i$. Note finally that  IVP one-- and multi--step  ODE integrators are not predictive, since they  use the right--hand side function data $f(t,y(t))$ of the ODE $y' = f(t,y(t))$ up to and including $t = t_{i+1}$ in order to evaluate $y'$'s antiderivative $y(t_{i+1})$.\\[1mm]
The onus of accuracy for classic general analytic continuation ODE solvers such as polygonal, Runge--Kutta, Heun, Prince--Dormand and embedding formulas lies within the chosen integration formula's accuracy for the solution data as their errors propagate and accumulate from each $t_i$ to $t_{i+1}$  and so forth in the computed solution $y(t)$ as time $t$ increases.\\[1mm]
All Zhang Neural Network based algorithms for time--varying matrix problem differ greatly in method, set--up, accuracy and speed from any  analytic continuation method that  computes antiderivative data for the given function $y' = f$ approximately, but without exponentially decreasing its  errors over time as ZNN does. 
}
 \end{minipage}
 \newpage.  
 
 \hspace*{4mm}\begin{minipage}{154mm}{
ZNN's fundamental difference lies in the stipulated (and actually computed and observed)  exponential error decay of the   solution.   ZNN solutions are rather impervious to noisy inputs, to occasional sensor failures and even to erroneous random data inputs (in a limited number of entries of the time--varying input matrix), see \cite{UZTVeigenLAA19} e.g., because an external disturbance will become invisible quickly in any solution whose error function decreases exponentially by its very design.\\[1mm]
A big unsolved challenge  for both  sensor driven matrix flow data acquisition and function based Zhang Neural Networks is their use of Taylor expansions, see formulas (10) through (15) in Section 2. Taylor expansions require smooth and multiply differentiable functions. But for   Zhang Neural Network algorithms and time--varying matrix problems  there is no guarantee of even once differentiable  inputs or solutions. Non--differentiable movements of robots  such as bounces are common events. Yet discretized ZNN is widely used with success for robot control.\\[0.5mm] 
 This dilemma seems to require a new understanding and theory of Taylor formulas for non--differentiable  time--varying functions.\\[1mm] 
The speed of error decay of ZNN based solutions hinges on a feasibly  chosen decay constant $\eta$ in Step 2 for the given sampling gap $\tau$ and on the truncation error order of the look-ahead convergent difference formula that is 
used. While analytical continuation is best used over small sub-intervals of increasing $t$ values inside $[a,b]$, we have run ZNN methods for seemingly near infinite time intervals $[0,t_f]$ with $t_f = 5, 10, 20, 100,  200, ..., 480,$ or $ 3600 $ seconds and up to eight hours and  have seen no deterioration in the error, or just the opposite, see Figure 1 above for the time--varying matrix square root problem, for example, whose relative error decreases continuously over a 6 minute span. Note further that time--varying matrix problems can be started  from almost any starting value $y_o$ in ZNN and the computed solution will very quickly lock onto the problem's proper solution.\\[1mm]
But ZNN cannot solve ODEs, nor can it solve static matrix/vector equations; at least we do not know how to. Here the term \emph{static matrix} refers to a matrix with constant real or complex entries.  On the other hand, path following IVP ODE based methods have been used to solve formulaic matrix equations successfully, see \cite{LM18, FUaccFOVLama20, FUDecompMatrFoV} e.g., but only implicit and explicit Runge--Kutta formulas  can  integrate sensor clocked data input while -- on the other hand -- suffering larger integration errors than adaptive integrators such as Prince--Dormand or ZNN.
} 
\end{minipage}\\[2.5mm]
Discretized Zeroing Neural Network methods and the quest for high order convergent and look-ahead finite difference formulas bring up many open  problems in numerical analysis:\\[1mm]
 Are there any look-ahead finite difference schemes with $ s < q$ in ${\cal A}_{q+s,q}$ and minimally  more rows than columns? Why or why not?\\[1mm]
For relatively low dimensions the rational numbers matrix  ${\cal A}_{q+s,q}$ in formula (16) above can easily be checked for  full rank when $1 \leq q \leq 6$. Is this true for all integers $q$? Has the rational ${\cal A}$  matrix in (16) ever been encountered anywhere else?\\[1mm]
 Every  look-ahead polynomial $p$ that we have constructed from any seed vector $y \in \RR^s$  with $s \geq q$ by our method has had precisely one root on the unit circle within $10^{-15}$ numerical accuracy.   This even holds for  non-convergent finite difference formula polynomials $p$ with some  roots outside  the unit disk. Do all Taylor expansion matrix ${\cal A}_{q+s,q}$ based polynomials have at least one root  on the unit circle in $\CC$?  What happens for polynomial finite difference formulas with all of their characteristic  roots  inside the open unit disk and none on the periphery or the outside? For the stability of multi--step ODE methods we must have $p(1) = \sum p_i = 0$, see \cite[p. 473]{EMU96} for example. \\[1mm]
For any low dimensional type {\tt j\_s} finite difference scheme there are apparently  dozens of convergent and look--ahead finite difference formulas for any fixed local truncation error order if $s \geq j > 1$.\\[0.5mm]
 What is the most advantageous such formula to use in  discretized ZNN  methods, speed-wise, convergency--wise?\\[0.5mm]
  What properties of a suitable formula improve  or hinder the ZNN computations for 
  time--varying matrix processes?\\[1mm]
 Intuitively we have preferred those convergent and look--ahead finite difference formulas whose characteristic polynomials have relatively  small second largest magnitude roots.\\[0.5mm]
  Is that correct and  a good strategy for discretized ZNN methods? See the graphs at the  end of \cite{FUaccFOVLama20} for examples.\\[0.5mm]
A list of further observations and open  problems for discretized ZNN based time--varying matrix eigen methods is included in \cite{UZTVeigenLAA19}.\\[1mm]
 The errors in ZNN's output come from   three sources, namely {\bf (a)} the rounding errors in solving the linear system in $(3i)$ or $(6i)$,  {\bf (b)} the truncation errors of the finite difference formula and the stepsize $\tau$ used, and {\bf (c)} from  the backward evaluation of the derivative $\dot A(t_k)$ in the right hand side expression of equation  (Az$i$)  or (Pqz$i$) above. How should one minimize or equilibrate their effects for  the best overall computed accuracy when using recurrence relations with high or low truncation error orders? 
High degree backward recursion formulas for derivatives are generally not very good.\\[-5mm]

\section{ Models and  Applications of Discretized ZNN Matrix Methods}

In this section we develop specific discretized ZNN algorithms for a number of selected time--varying matrix problems.  Moreover we introduce new matrix techniques to transform  matrix  models  into Matlab code and we  point to problem specific references.\\[2mm]
Our first example deals with the  oldest human  matrix problem, namely solving linear equations $A_{n,n} x = b_n$. This model goes back well over 6,000 years to Babylon and Sumer on cuneiform tablets that describe Gaussian elimination techniques and solve static linear equations $Ax=b$ for small dimensions $n$.
\\

{\bf (I) \ Time--varying Linear Equations and Discretized ZNN :}\\[1mm]
For simplicity, we consider matrix flows $A(t)_{n,n} \in \CC_{n,n}$ all of whose matrices are invertible on a time interval $t_o \leq t \leq t_f \subset \RR$. Our chosen model equation is  \circledzero $A_{n,n} (t)x(t) = b(t) \in \CC^n$ for the unknown solution vector $x(t)$. The  error function is  \circledone  $e(t) = A(t)x(t) - b(t)$ and the error differential equation is  
$$ \circledtwo  \ \dot e(t) \ = \dot A(t)x(t) + A(t)\dot x(t) - \dot b(t) \stackrel{(*)}{=} -\et~ A(t)x(t) + \et~ b(t) = -\et~ e(t) \ .
$$
Solving the inner equation $(*)$ in \circledtwo first for $A(x)\dot x(t)$  and then for $\dot x(t)$ we obtain  the the following two differential equations (DEs,) see also  \cite[(4.4)]{ZW01}
$$  A(t) \dot x(t) = - \dot A(t) x(t) + \dot b(t) -\et~ A(t)x(t) +\et~  b(t)\ $$
and
$$ \circledthree \ \dot x(t) = A(t)^{-1} \ (  - \dot A(t) x(t) + \dot b(t) + \et~  b(t) ) -\et~ x(t)\ .$$
 To simplify matters we use the simple 5--IFD formula (11$i$) of Section 2 again in discretized mode with $A_k = A(t_k), x_k = x(t_k)$ and $b_k = b(t_k)$ to obtain 

$$ \circledfour  \ \ \dot x_k =                                                                                                                                                                                                                                                                                                                                                                 \dfrac{8x_{k+1}+x_k -6x_{k-1} - 5 x_{k-2} + 2 x_{k-3}}{18 \tau}  \in \CC^{n} \ .
$$
Equating derivatives at time $t_k$ yields 
$$ \circledfive \  \ 18\tau \cdot \dot x_k =  8x_{k+1}+x_k -6x_{k-1} - 5 x_{k-2} + 2 x_{k-3} \stackrel{(*)}{=} 18\tau \cdot \left(  A^{-1}_k \ (  - \dot A_k x_k + \dot b_k + \et~ \ b_k ) -\et~ \ x_k \right) \ .\ 
$$
Then the inner equation $(*)$ of \circledfive  gives us the predictive convergent and look-ahead ZNN formula  
$$ \circledsix  \  \  x_{k+1} = \dfrac{9}{4} \tau \cdot \left( A^{-1}_k \ (  - \dot A_k x_k + \dot b_k + \et~  b_k ) -\et~  x_k  \right)- \dfrac{1}{8} x_k + \dfrac{3}{4} x_{k-1} + \dfrac{5}{8} x_{k-2} - \dfrac{1}{4} x_{k-3} \in \CC^n \  .
$$
Since equation \circledsix  involves the matrix inverse $A^{-1}_k$ at each time step $t_k$, we propose  two different  Matlab codes to solve time--varying linear equations for invertible matrix flows $A(t)$. The code  {\tt tvLinEquatexinv.m} in \cite{FUZNNSurvey} uses  Matlab's matrix inversion method {\tt inv.m} explicitly at each time step $t_k$ as  \circledsix requires, while our second code {\tt tvLinEquat.m}  in \cite{FUZNNSurvey} uses two separate discretized ZNN formulas. 
The {\tt tvLinEquat.m}  code solves the time--varying linear equation with the help of one  ZNN method that computes the inverse of each $A(t_k)$ iteratively as detailed next in Example {\bf (II)} below and another ZNN method that that solves equation \circledsix by using these two independent  and interwoven discretized ZNN iterations.\\[1mm] 
Both methods run equally fast. The first  with its  explicit matrix inversion  is a little  more accurate since  Matlab's {\tt inv} computes small dimensioned matrix inverses  near perfectly with $10^{-16}$ relative errors while   ZNN based time--varying matrix inversions give us  slightly larger errors, losing 1 or 2 accurate trailing digits. This is most noticeable if we use low truncation error order look--ahead finite difference formulas for ZNN and relatively large sampling gaps $\tau$. 
There are dozens of references when googling 'ZNN for time--varying linear equations', see also \cite{XLNSG19} or \cite{ZYYHXH19}. \\[2mm]
{\bf (II) \ Time--varying Matrix Inverses via ZNN :}\\[1mm]
Assuming  again that all matrices of a given time--varying matrix flow $A(t)_{n,n} \in \CC_{n,n}$ are invertible on a given time interval $t_o \leq t \leq t_f \subset \RR$,  we construct a discretized ZNN method that  finds the inverse $X(t)$  of each $A(t)$  predictively from previous data so that $A(t)X(t) = I_n$, i.e.,  \circledzero $A(t) = X(t)^{-1}$  is our model here. This gives rise to the error function \circledone $E(t) = A(t) - X(t)^{-1} \ (= O_{n,n} \text{ ideally})$ and the associated error differential equation 
$$ \circledtwo  \ \dot E(t) = \dot A(t) - \dot X(t)^{-1} \ .$$
Since $X(t)X(t)^{-1} = I_n$ is constant for all $t$, $d(X(t)X(t)^{-1})/dt = O_{n,n}$. And the product rule gives us the following relation  for the derivative of time--varying matrix inverses  
$$ O_{n,n} =\dfrac{d(X(t)X(t)^{-1})}{dt} = \dot X(t) X(t)^{-1} + X(t) \dot X(t)^{-1} \ ,$$ 
 and thus $\dot X(t)^{-1} = -X(t)^{-1} \dot X(t) X(t)^{-1}$\ . 
Plugging this derivative  formula   into \circledtwo establishes 
$$   \dot E(t) = \dot A(t) - \dot X(t)^{-1} = \dot A(t) + X(t)^{-1}\dot X(t) X(t)^{-1} \stackrel{(*)}{=} -\et~ A(t) + \et~ X(t)^{-1} = -\et~ E(t). $$
Multiplying  the inner equation $(*)$ above by $X(t)$ from the left on both, the left-- and right--hand sides and then solving for $\dot X(t)$ yields
$$ \circledthree \ \dot X(t) = - X(t) \dot A(t) X(t) - \et~ X(t)A(t)X(t) + \et~ X(t) = -X(t)( ( \dot A(t) + \et~ A(t) X(t) - \et~ I_n)\ .$$
If -- for simplicity --  we choose the same 5--IFD look--ahead and convergent formula as was chosen on line (4$i$) for Step 4 of the ZNN eigen--data method in Section 1, then we obtain the analogous equation to (12$i$) here with $(P\backslash q)$ replaced by the 
right--hand side of equation \circledthree $\!\!$.  Instead of  the eigen--data iterates $z_j$ in (12$i$) we use the inverse matrix iterates $X_j = X(t_j)$ here for $j = k-3, ...,k+1$ and obtain
$$ \circledfive\ 18\tau \cdot \dot X_k = 8 X_{k+1} + X_k- 6 X_{k-1} - 5 X_{k-2} + 2 X_{k-3} \stackrel{(*)}{=} - 18 \tau \cdot X_k( ( \dot A(t_k) + \et~ A(t_k) X_k - \et~ I_n)\ .$$
Solving $(*)$ in \circledfive for $X_{k+1}$ supplies the complete ZNN recursion formula that finishes Step 6 of the predictive discretized ZNN algorithm development for 
time--varying matrix inverses.
$$ \circledsix X_{k+1} = - \dfrac{9}{4} \tau~  X_k( ( \dot A(t_k) + \et~ A(t_k)X_k - \et~ I_n )- \dfrac{1}{8} X_k + \dfrac{3}{4} X_{k-1} + \dfrac{5}{8} X_{k-2} - \dfrac{1}{4} X_{k-3} \ \in \CC_{n,n} $$
This look--ahead iteration is based on the convergent 5--IFD  formula  of type {\tt j\_s = 2\_3} with local truncation error order $O(\tau^4)$. The formula \circledsix requires two matrix multiplications, two matrix additions, one backward approximation of the derivative of $A(t_k)$ and a short recursion with $X_j$ at each time step. \\
The error function differential equation \circledthree is akin to the Getz and Marsden dynamic system (without the discretized ZNN $\eta$ decay terms)  for time--varying matrix inversions, see \cite{GM97} and \cite{GZ12}. Simulink circuit diagrams for this model and time--varying matrix inversions are available in \cite[p. 97]{ZGbook15}.\\
A Matlab code for the time--varying matrix inversion problem is available in \cite{FUZNNSurvey} as {\tt  tvMatrixInverse.m}. A different model is used in \cite{ZLL11} and several others are described in \cite[chapters 9, 12]{ZGbook15}.\\[-2mm]
\begin{center}
\includegraphics[width=128mm]{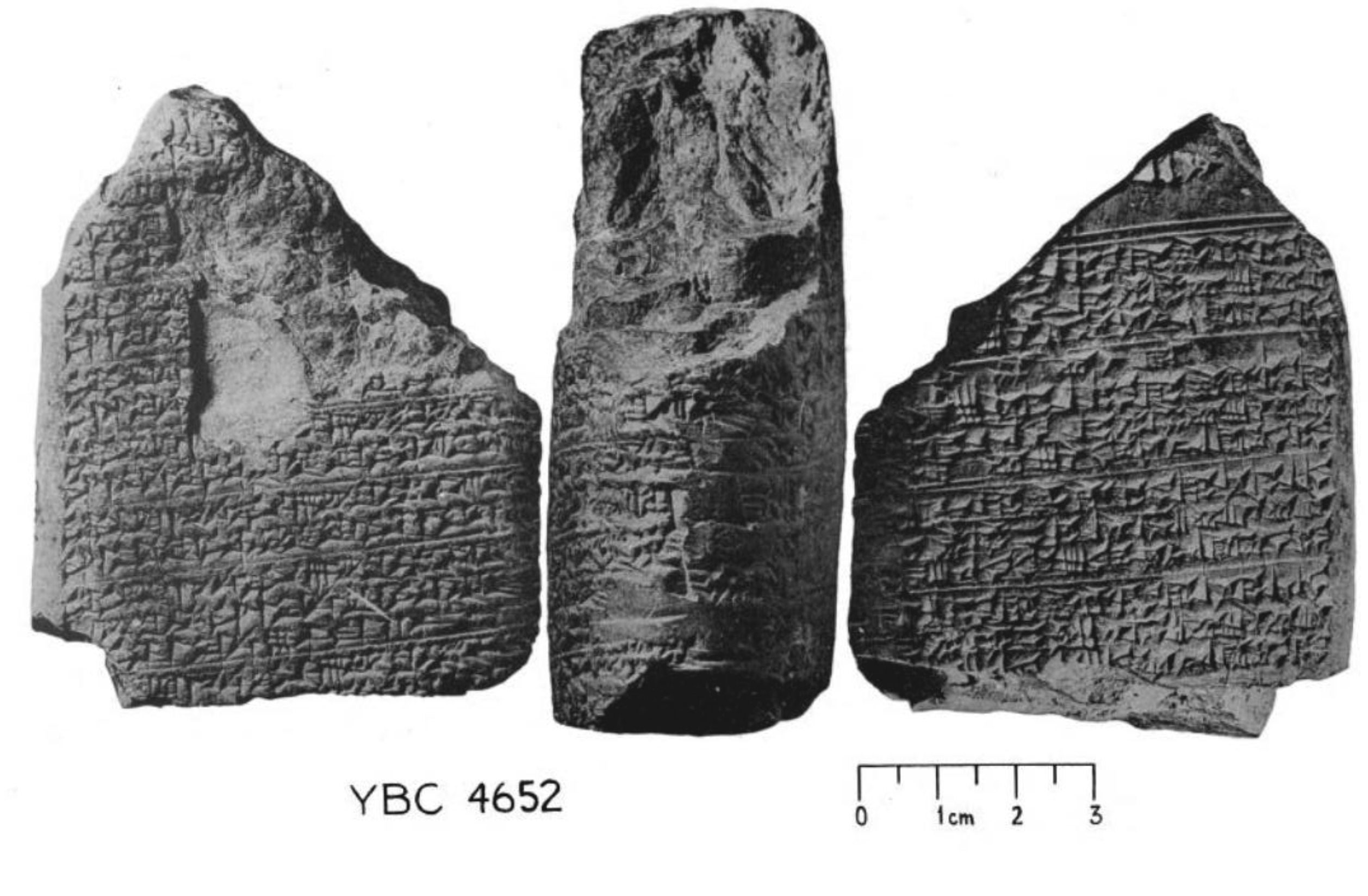}  \\[-6mm]
Figure 2  : Cuneiform tablet (from Yale) with Babylonian methods for solving a system of two  linear equations.\\
{[} Search for  \ CuneiformYBC4652 \ first to then learn more on YBC 04652 from the Cuniform Digital Library Initiative in Berlin at \ {\url {https://cdli.mpwg-berlin.mpg.de}} \  and from other articles on YBC 4652.\ ] \end{center}

\vspace*{-1mm}
The image above shows how ubiquitous matrices and matrix computations have been  across the eons, cultures and languages of humankind on our globe. It took years to learn and translate  cuneiform symbols and  to realize that 'Gaussian elimination' was already well  understood and used then. And it is still central for matrix computations today. \\[1.5mm]
{\bf Remark 2 : } (a) Example {\bf (II)} reminds us that  the numerics for time--varying matrix problems may differ greatly from our static matrix numerical approaches for matrices  with constant entries. Time--varying matrix problems are governed by different concepts and follow different best practices.\\
 For static matrices $A_{n,n}$ we are always conscious of and we remind our students never to compute the inverse $A^{-1}$ in order to solve a linear equation $Ax = b$ because this is an expensive proposition and  rightly shunned.\\
 But for time--varying matrix flows $A(t)_{n,n}$ it seems impossible to solve 
time--varying linear systems $A(t)x(t) = b(t)$ predictively without explicit matrix inversions as was  explained in part (I) above. For time--varying linear equations, ZNN methods allow us  to compute  time--varying matrix inverses and solve time--varying linear equations  in real time, accurately and  predictively. What is shunned for  static matrix problems may work  well for the time--varying matrix variant and vice versa.\\[1mm]
(b) For each of the example problems in this section our annotated rudimentary ZNN based Matlab codes are stored in \cite{FUZNNSurvey}. For high truncation error order look-ahead convergent finite difference formulas such as {\tt 4\_5} these codes achieve 12 to 16 correct leading digits predictively for each entry of the desired solution matrix or vector and they do so uniformly for all parameter values of $t$ after the initial  exponential error reduction has achieved this accuracy.\\[1mm]
(c) {\bf General warning} : Our published codes in \cite{FUZNNSurvey} may not apply to all uses and may give incorrect results for some inputs and in some applications. These codes are built here explicitly only for educational purposes. Their complication level advances in numerical complexity as we go on  and our explicit codes try to help and show  users how to  solve time--varying applied and theoretical matrix problems with discretized ZNN methods. We strongly advise users to search the theoretical literature on their specific problem and to test their own applied ZNN codes rigorously and extensively before applying them in the field. This scrutiny will ensure that theory or code exceptions do not cause unintended consequences or accidents in the field.\\[3mm] 
{\bf (III) \ Pseudo--inverses of Time--varying Non--square Matrices with Full Rank and without :}\\[1mm]
 Here we first  look at rectangular matrix flows $A(t)_{m,n} \in \CC_{m,n}$ with  $m \neq n$ that have  uniform full  rank($A(t)) =$ min$(m,n)$ for all $t_o \leq t \leq t_f \subset \RR$.\\
 Every matrix $A_{m,n}$ with $ m = n$ or $m \neq n$ has two kinds of nullspaces or kernels 
 $$N(A)_r = \{x \in \CC^n \mid Ax = 0 \in \CC^m\} \ \ \ \text{and} \ \ \ 
 N(A)_\ell = \{x \in \CC^m \mid xA = 0 \in \CC^n\} \ ,$$ 
 called  $A$'s  right and left nullspace, respectively. If $m>n$ and $A_{m,n}$ has full rank $n$, then $A$'s right kernel is $\{0\} \subset \CC^n$ and the linear system $Ax=b \in \CC^m$ cannot be solved for every $b \in \CC^m$  since the number the columns of $A_{m,n}$ is less than required for spanning all of $\RR^m$ . If $m < n$ and $A_{m,n}$ has full rank $m$, then $A$'s left kernel is $\{0\} \subset \CC^m$ and similarly not all equations $xA=b \in \CC^n$ are solvable with $x \in \CC^m$. Hence we need to abandon the notion of  matrix inversion for rectangular non-square matrices and resort to pseudo-inverses instead.\\
  There are two kinds of  pseudo-inverses of $A_{m,n}$, too, depending on whether $m >n$ or $m<n$. They are always denoted by $A^+$ and always have size $n$ by $m$ if $A$ is $m$ by $n$.  If $m>n$ and $A_{m,n}$ has full rank $n$, then $A^+ = (A^*A)^{-1}A^* \in \CC_{n,m}$ is called the  left pseudo-inverse because $A^+A = I_n$. For $m<n$ the right pseudo-inverse of $A_{m,n}$ with full rank $m$  is  $A^+ = A^*(AA^*)^{-1}  \in \CC_{n,m}$ with $AA^+ = I_m$. \\[0.3mm]
    In either case $A^+$ solves a minimization problem, i.e.,  
    $$\min_{x\in \CC^n} ||Ax - b||_2 = ||A^+ b||_2 \geq 0 \ \ \text{ for } \ m>n \ \text{ and } \ \min_{x \in \CC^m} ||xA - b||_2 = ||bA^+||_2 \geq 0 \ \  \text{ for } \  m<n .
    $$
   
     Thus the pseudo--inverse of a full rank rectangular matrix $A_{m,n}$ with $m \neq n$ solves the least squares problem for sets of linear equations whose  system matrices $A$ have nontrivial left or right kernels, respectively. It is easy to verify that  $(A^+)^+ = A$ in either case, see e.g. \cite[section 4.8.5]{SB02}. Thus the pseudo--inverse  $A^+$ acts similarly to the matrix inverse $A^{-1}$ when  $A_{n,n}$ is invertible and $m=n$. Hence its name.\\[1mm]
First we want to find the  pseudo--inverse $X(t)$ of a full rank rectangular matrix flow $A(t)_{m,n}$ with $m < n$.
Since $X(t)^+ = A(t)$  we can try to  use the dynamical system of Getz and Marsden \cite{GM97} again  and start with \ \circledzero  $A(t)=X(t)^+$ as our model equation.\\[2mm]
{\bf \emph{(a) \ The right pseudo--inverse model \ \circledzero $X(t)  = A(t)^+$ for  matrix flows $A(t)_{m,n}$ of full rank $m$ when $m<n$} :}\\[1mm]
The exponential decay stipulation for our model's error function  \circledone $E(t) = A(t) - X(t)^+$ makes 
$$ \circledtwo \ \ \dot E(t) = \dot A(t) - \dot X(t)^+ \stackrel{(*)}{=} -\et~  A(t) + \et~  X(t)^+ = -\et\  E(t)  \ .$$ 
Since $A(t)X(t) = I_m$ for all $t$ and $A(t) = X(t)^+$ we have  
$$O_{m,m} = d(I_m)/dt = d(A(t)X(t))/dt = d(X(t)^+X(t))/dt = \dot X(t)^+ X(t) + X(t)^+ \dot X(t) \ .$$
Thus $\dot X(t)^+  X(t) = - X(t)^+ \dot X(t)$ or $ \dot X(t)^+ = - X(t)^+ \dot X(t) X(t)^+$ after multiplying equation $(*)$ in \circledtwo  by $X(t)^+$ on the right.
Updating  equation \circledtwo establishes
$$ \dot A(t) + X(t)^+\dot X(t) X(t)^+= - \et~  A(t) + \et~  X(t)^+  \ .
$$
Multiplying both sides on the left and right by $X(t)$ then yields
$$  X(t)\dot A(t) X(t)+ X(t)X(t)^+\dot X(t) X(t)^+ X(t)= - \et~  X(t) A(t) X(t)+ \et~  X(t) X(t)^+  X(t) \ .
$$
Since $X(t)^+ X(t) = I_n$ we obtain after  reordering that 
$$  X(t)X(t)^+\dot X(t) = -  X(t)\dot A(t) X(t) - \et~  X(t) A(t) X(t)+ \et~  X(t) X(t)^+  X(t) \ .
$$
But $X X(t)^+ $ has size $n$ by $n$ and rank $m<n$. Therefore it is not invertible. And   thus we cannot cancel the encumbering left factors for $\dot X(t)$ above  and solve the  equation  for  $ \dot X(t)$ as would be needed for Step 3. And a valid ZNN formula cannot be  obtained from our first simple model $A(t) = X(t)^+$.\\[1mm]
This example contains a warning not to give up if one model does not work for a 
time--varying matrix problem.\\[1mm]
Next we  try another model equation for the right pseudo--inverse  $X(t)$ of a full rank matrix flow $A(t)_{m,n}$ with $m<n$. Using the definition of $A(t)^+ = X(t) = A(t)^*(A(t)A(t)^*)^{-1}$ we start from the revised model \circledzero $X(t) A(t)A(t)^* = A(t)^*$.
With the error function \circledone $E = XAA^* - A^*$ we  obtain (leaving out all time dependencies of $t$  for better readability) 
$$ \circledtwo \ \ \ \dot E = \dot XAA^* + X \dot AA^* +XA\dot A -\dot A^* \stackrel{(*)}{=} -\et~ XAA^* + \et~ A^* = - \et~ E \ .
$$
Separating the term with the unknown derivative $\dot X$ on the left of $(*)$ in \circledtwo$\!\!$, this becomes 
$$
 \dot X \ AA^* = -X\left( (\dot A + \et~ A)A^* + A \dot A^*\right) + \dot A^* + \et~ A^* \ .
$$
Here the matrix product $A(t)A(t)^*$on the left--hand side  is of size $m$ by $m$ and  has rank $m$ for all $t$ since $A(t)$ does. Thus  we have found an explicit expression for $\dot X(t)$, namely
 $$
 \circledthree \ \ \ \dot X  = \left(-X\left( (\dot A + \et~ A)A^* + A \dot A^*\right) + \dot A^* + \et~ A^*\right) (AA^*)^{-1} \ .
$$
The steps \circledfour$\!$, \circledfive and \circledsix now follow as before. The Matlab ZNN based discretized code for right pseudo-inverses is {\tt tvRightPseudInv.m} in \cite{FUZNNSurvey}. Our code finds  right pseudo--inverses of time--varying full rank matrices $A(t)_{m,n}$ predictively with an entry accuracy of 14 to 16 leading digits in every position of $A^+(t) = X(t)$ when compared with the pseudo-inverse defining formula.  In the code we use the {\tt 4\_5} look-ahead convergent finite difference formula from \cite{FUfindiff19} with the sampling gap $\tau = 0.0002$.\\[1mm]
 Similar numerical results are obtained for left pseudo--inverses $A(t)^+$ for time--varying matrix flows $A(t)_{m,n}$ with $m>n$.\\[2mm]
{\bf \emph{(b) \ The left pseudo--inverse $X(t)  = A(t)^+$ for matrix flows 
$A(t)_{m,n}$ of full rank $n$ when $m>n$} :}\\[1mm]
Our starting model now is \circledzero $ A^+ = X_{n,m} = (A^*A)^{-1}A^*$ and the error function \circledone  $ E = (A^*A)X - A^*$ then leads to
$$ \circledtwo  \ \ \ \dot E = \dot A^*AX + A^*\dot AX + A^* A \dot X - \dot A^* = - \et~ A^* A X + \et~ A^* = -\et~ E \ .
$$
Solving \circledtwo for $\dot X$ similarly as before yields
$$ \circledthree \ \ \  \dot X = (A^*A)^{-1} \left(-\left((\dot A^* + \et~ A^*) + A^* \dot A\right) X + \dot A^* + \et~ A^*\right) \ .
$$
Then follow the steps from subpart (a) and develop a Matlab ZNN code for left 
pseudo--inverses with a truncation error order finite difference formula of your own choice.\\[1mm]
{\bf \emph{(c) \ The  pseudo--inverse \circledzero $X(t)  = A(t)^+$ for matrix flows $A(t)_{m,n}$ with variable rank$\mathbf {(A(t)) \ \leq \ \min(m,n)}$} :}\\[1mm]
As before with the unknown pseudo--inverse $X(t)_{n,m}$ for a possibly rank deficient matrix flow $A(t) \in \CC_{m,n}$, we use the error function \circledone $E(t) = A(t) - X(t)^+$ and the error function DE
$$ \circledtwo \ \ \ \dot E(t) = \dot A(t) - \dot X(t)^+ \stackrel{(*)}{=} -\et~  A(t) + \et~  X(t)^+ = -\et~  E(t)  \ .$$ 
For  matrix flows $A(t)$ with rank deficiencies  the derivative of $X^+$, however,  becomes more complicated with additional terms, see \cite[Eq. 4.12]{GP73} :
\begin{equation}\label{longpinv} \dot X^+ = -X^+\dot XX^+ + X^+X^{+^*}\dot X^* (I_n-XX^+) + (I_m-X^+X)\dot X^* X^{+^*}X^+ 
\end{equation}
where previously for full rank matrix flows $A(t)$, only the first term above was  needed to express  $\dot X^+$. Plugging the long expression in (\ref{longpinv}) for $\dot X^+$ into  the inner equation (*) of \circledtwo we obtain 
$$ \circledthree \ \ \dot A + X^+\dot XX^+ - X^+X^{+^*}\dot X^* (I_n-XX^+) - (I_m-X^+X)\dot X^* X^{+^*}X^+ =  -\et~  A(t) + \et~  X(t)^+ $$
which needs to be solved for $\dot X$. Unfortunately $\dot X$ appears once  on the left in the second  term and twice as $\dot X^*$ in the third and fourth term of \circledthree above. Maybe another start--up error function can give better results, but it seems that the general rank pseudo--inverse problem cannot be easily solved via the ZNN process, unless we learn to work with Kronecker matrix products. Kronecker products   will be used in  subparts ({\bf VII}), ({\bf IX}) and ({\bf VII start--up}) below.\\[1mm]
The Matlab code {\tt  tvLeftPseudInv.m}   for ZNN look--ahead left pseudo-inverses of full rank  time--varying matrix flows is available in  \cite{FUZNNSurvey}. The right 
pseudo--inverse code for full rank matrix flows is similar. Recent work on pseudo--inverses has appeared in \cite{SWM19} and \cite[chapters 8,9]{ZGbook15}.\\[-2mm]

{\bf (IV) \ Least Squares, Pseudo--inverses and ZNN :}\\[2mm]
Linear systems of time--varying equations $A(t)x(t) = b(t)$ can be unsolvable or solvable with unique or multiple solutions and pseudo--inverses can help us.\\
 If the matrix flow $A(t)_{m,n}$ admits a left pseudo--inverse $A(t)^+_{n,m}$ then  
 $$A(t)^+A(t)x(t) = A(t)^+b(t) \ \ \text{ and } \ \  x(t)=(A(t)^+A(t))^{-1} A^+(t)b(t) \ \ \text{  or } \ \ x(t) = A(t)^+b(t) \ .$$ 
 Thus $A(t)^+b(t)$  solves the linear system at each time $t$ and $x(t) = A(t)^+b(t)$  is the solution with minimal Euclidean norm $||x(t)||_2$ since according to (\ref{longpinv}) all other time--varying solutions have the form $$u(t) = A(t)^+b(t) + (I_n - A(t)^+A(t))w(t) \ \ \text{ for any } \ \ w(t) \in \CC^m \ .$$  Here $||A(t)x(t) - b(t)||_2 = 0$ holds precisely when $b(t)$ lies in the span of the columns of $A(t)$ and the linear system is uniquely solvable. Otherwise $\min_x(||A(t)x(t) - b(t)||_2) > 0$.\\[1mm] 
Right pseudo--inverses $A(t)^+$ play the same role for linear-systems $y(t)A(t) = c(t)$. In fact $$y(t)A(t)A(t)^* = c(t)A(t)^* \ \ \text{  and } \ \ y(t) = c(t)A(t)^*(A(t)A(t)^*)^{-1} = c(t)A(t)^+ \ . $$
Here $c(t)A(t)^+$ solves the left-sided linear system $y(t)A(t) = c(t)$ with minimal Euclidean norm.\\[1mm]
 In this subsection  we will only work on time--varying linear equations of the form  \circledzero $A(t)_{m,n} x(t) = b(t) \in \CC_m$ for $m > n$ with rank$(A(t)) = n$ for all $t$. Then the left pseudo--inverse of $A(t) $ is $A(t)^+=  (A(t)^*A(t))^{-1}A^*$. The associated error function is \circledone $e(t) = A(t)_{m,n} x(t) - b(t)$. Stipulated exponential error decay defines the error function DE\\[-4mm]
$$
\circledtwo \ \ \ \dot e = \dot A x + A \dot x - \dot b \stackrel{(*)}{=} -\et~ Ax + \et~  b = - \et~ e
$$
where we have again left off the time parameter $t$ for clarity and simplicity. Solving $(*)$ in \circledtwo for $\dot x(t_k)$ gives us
$$
\circledthree \ \ \ \dot x_k = (A_k^*A_k)^{-1} A^* \left( -(\dot A_k + \et~ A_k)x_k + \dot b_k - \et~ b_k\right) \ .
$$
Here the subscripts $\dots_{k}$ remind us that we are describing the discretized version of our Matlab codes where $b_k$ for example stands for $b(t_k)$ and so forth for $A_k, x_k, ... $. The Matlab code for the discretized  ZNN look-ahead  method for time--varying linear equations  least squares problems for full rank  matrix flows $A(t)_{m,n}$ with $m > n$ is  {\tt  tvPseudInvLinEquat.m} in \cite{FUZNNSurvey}. We advise readers to develop a similar code for full rank  matrix flows $A(t)_{m,n}$ with $m < n$ independently.\\[1mm]
 The survey article  \cite{LSX20} describes nine discretized ZNN methods for  time--varying different matrix optimization problems such as least squares and constrained optimizations that we treat in  subsection (V) below.  \\[-2mm]

{\bf (V) \ Linearly Equality Constrained Nonlinear Optimization for Time--varying  Matrix Flows :}\\[2mm]
ZNN can be used to solve parametric nonlinear programs (parameteric NLPs). However, ensuring that the solution path exists and that the extremum is  isolated for all $t$ is nontrivial, requiring careful tracking of the active set and avoiding various degeneracies. Simpler sub-classes of the general problem can readily be solved without the extra machinery. For example, optimization problems with $f(x(t), t)$ nonlinear but only linear equality constraints and no inequality constraints, such as 
\begin{equation} \text{\circledzero} \ \ \text{ find } \ \min f(x(t), t) \in \RR\ \ \ \text{subject to} \ \ \  A(t) x(t) = b(t)  \label{eqn:nlplinconstr}  \end{equation}

for $f :~ (\RR^n \times \RR) \to \RR$, $x(t) \in \RR^n$, $\la(t) \in \RR^m$, and $A(t) \in \RR_{m,n}$ and $b(t) \in \RR^n$ acting as linear equality constraints $A(t)x(t) = b(t)$, we can write out the corresponding Lagrangian (choosing the + sign convention),
$$ {\cal L}(x(t),\la(t)) = f(x(t), t) + \la(t)  (A(t)x(t) - b(t)) \ : \ (\RR^{n} \times \RR^{m}) \to \RR \ .$$

\vspace*{-1mm}
A necessary condition for ZNN to be successful is that a solution $x^*(t)$ must exist and that it is an isolated extremum for all $t$. Clearly the Lagrange multipliers, $\la^*(t)$, must also be unique for all $t$. Furthermore, due to Step 3 in the ZNN setup, both solution $x^*(t)$ and $\la^*(t)$ must have a continuous first derivatives. Therefore we further suppose that $f(x(t), t)$ is twice continuously differentiable with respect to $x$ and $t$, and that $A(t)$ and $b(t)$ are twice continuously differentiable with respect to $t$.\\[2mm]
To ensure that the $\la^*(t)$ are not only inside  a bounded set but are also unique \cite{Gau77}, LICQ must hold, i.e.,  $A(t)$ must have  full rank for all $t$. If (\ref{eqn:nlplinconstr}) is considered for static entries, i.e., for some fixed $t=t_0$, and we suppose that $x^*(t)$ is a local solution then there exists a Lagrange multiplier vector $\la^*(t)$ such that the first order necessary conditions \cite{BNO03, NW06}, or `KKT conditions', are satisfied:
\begin{subequations}
\begin{align}
\label{lagreq}\nabla_{x} {\cal L}(x^*(t),\la^*(t)) = \nabla_{x} f(x^*(t)) + (\la^*(t))^T A(t) &= 0 \in \RR^n,\\
\label{eqn:lincnstr} A(t) x^*(t) - b(t) &=0  \in \RR^m .
\end{align}
\end{subequations}

\vspace*{-2.5mm}
Here $\nabla_{x}$ denotes the gradient, in column vector form. A second order sufficient condition must also be imposed here that $y^T \nabla_{xx}^2  {\cal L}(x^*(t),\la^*(t)) y > 0$ for all $y \ne 0$ with $Ay = 0$. This ensures the desired curvature of the projection of $f(\cdot)$ onto the constraints. Theorem 2.1 from \cite{Fi76} then establishes that the solution exists, is an isolated minimum, and is continuously differentiable. For the reader's reference, further development of these ideas can be seen in, e.g.\ \cite{RD95}, which considers the general parametric NLP under the Mangasarian--Fromovitz Constraint Qualification (MFCQ).\\[2mm]
The system of equations (\ref{lagreq}, \ref{eqn:lincnstr}) is the starting point for the setup of ZNN and this optimization problem. For notational convenience, the $x^*(t)$ notation will be simplified to $x(t)$ since it is clear what is meant. We want to solve for time--varying $x(t)$ and $\la(t)$ in the predictive discretized ZNN fashion.  I.e., we want to find $x(t_{k+1})$ and $\la(t_{k+1})$ from earlier data for times $t_j$ with $ j = k, k-1, ...$ accurately in real time. First we define $y(t) := {\tt [x(t);\la(t) ]} \in \RR^{n+m}$ in Matlab column vector notation and use (\ref{lagreq}) and (\ref{eqn:lincnstr})}to define our error function as \\[-7mm]

$$ \circledone \ \ h(y(t),t) := \bp \nabla_x {\cal L}(x(t),\la(t)) \\ 
A(t) x(t) - b(t) \ep = \bp \nabla_x f(x(t))  + (\la(t))^T A(t) \\ A(t)x(t) - b(t)\ep
= \bp h_1(t)\\ \vdots\\ \vdots\\  h_{n+m}(t) \ep.  \ 
$$

\vspace*{-1mm}
To find the derivative  $\dot y(t)$ of $y(t)$ we use the multi--variable chain rule which establishes  the  derivative of $h(y(t))$ as\\[-2mm]
\begin{equation}
    \dot{h}(x(t), t) = \bp \left(\nabla_{xx}^2 f(x(t), t) \right) \dot{x} & A(t)^T \\ A(t) & 0 \ep  \dot{y} + \bp \nabla_x f_t(x(t), t) + \dot{A}(t)^T \la \\ \dot{A}(t) x(t) - \dot{b}(t) \ep , \label{eqn:stdpfmtrxs}
\end{equation}

where the $\nabla_{xx}^2$ denotes the Hessian. The expression in (\ref{eqn:stdpfmtrxs}) is a slight simplification of a more general formulation often used in parametric NLPs when the active set is fixed; usually, it is used for numerical continuation directly by setting it equal to zero and integrating, but here we are interested in using it with ZNN.\\[1mm]
In our restricted case, we could use an equivalent formulation and  suppose that the equality constraints are arising from the Lagrangian by taking gradients with respect to both $x$ and $\lambda$ as follows\\[-6mm]

$$ \dot h(y(t),t) =  J(h(y(t),t)) \ \dot y(t) + h_t(y(t), t)\ . $$

\newpage

Here \\[-4mm]
$$  J(h(y(t),t)) =  \bp \dfrac{\partial h_1(t)}{\partial x_1} & \cdots & \dfrac{\partial h_1(t)}{\partial \lambda_m}\\  
\vdots & & \vdots\\
\dfrac{\partial h_{n+m}(t)}{\partial x_1} &\cdots & \dfrac{\partial h_{n+m}(t)}{\partial \lambda_m} \ep_{n+m,n+m}   \ \ \ \text{ and } \ \ \ h_t(y(t)) = \bp \dfrac{\partial h_1(y(t))}{\partial t}\\ \vdots \\ 
\dfrac{\partial h_{n+m}(y(t))}{\partial t}\ep_{n+m}
$$

are the Jacobian matrix $J$ of $h(y(t),t)$ taken with respect to the location vector $x(t) = (x_1(t), \dots, x_n(t))$ and the Lagrange multiplier vector $(\lambda_1(t), \dots, \la_m(t))$, and the time derivative of $h(y(t),t)$, respectively. This formulation, however, does not apply when we encounter  inequality constraints that are non-differentiable  on a set of measure zero or further difficulties with active and inactive constraints, etc. For otherwise 'regular' linear time-varying optimization problems we simply start from the standard Lagrangian `Ansatz' : \\[-6mm] 

$$ \text{\circledtwo} \ \dot h(y(t),t) = - \et~ h(y(t),t)$$

\vspace*{-0mm}
   which will lead us exponentially fast to  the optimal solution $y(t)$ for $t_o \leq t \leq t_f$.  Solving for $\dot y(t)$  gives us\\[-6mm]
   
$$ \circledthree   \ \dot y(t) = - J(h(y(t)),t)^{-1} \left( \et~ h(y(t),t) + \dot h_t(y(t),t )\right) \ .
 $$
 
 \vspace*{-0mm}
Using the 5--IFD look-ahead finite difference formula once again, this time  for $\dot y_k$ with discretized data $y_k = y(t_k)$, we obtain the following solution--derivative free equation for the iterates $y_j$ with $j \leq k$  by equating the two expressions for $\dot y_k$ in the 5-IFD discretization formula and in \circledthree as follows : \\[-2mm]
$$ \circledfive\ 18\tau \cdot \dot y_k = 8 y_{k+1} + y_k- 6 y_{k-1} - 5 y_{k-2} + 2 y_{k-3} \stackrel{(*)}{=} - 18 \tau \cdot \ J(h(y_k))^{-1} \left( \et~ h(y_k)) + \dot h_t(y_k) \right)= 18\tau\cdot \dot y_k  .$$

\vspace*{-1mm}
Solving $(*)$ in \circledfive for $y_{k+1}$ supplies the complete discretized ZNN recursion formula that finishes Step 6 of the predictive discretized ZNN algorithm development for 
time--varying constrained non--linear optimizations via  Lagrange multipliers:\\[-2mm]
$$ \circledsix \ y_{k+1} = - \dfrac{9}{4} \tau \cdot J(h(y_k))^{-1} \left( \et\ h(y_k)) + \dot h_t(y_k) \right) - \dfrac{1}{8} y_k + \dfrac{3}{4} y_{k-1} + \dfrac{5}{8} y_{k-2} - \dfrac{1}{4} y_{k-3} \ \in \CC^{n+m} \ .$$
The Lagrange based  optimization algorithm for multivariate functions and constraints is coded for one specific example with $m = 1$  and $n = 2$ in {\tt tvLagrangeOptim2.m}, see \cite{FUZNNSurvey}. For this specific example  the optimal solution is known. The code can be modified for optimization problems with more than $n=2$ variables and for more than $m=1$  constraint functions. 
Our code is modular and accepts all look-ahead convergent finite difference formulas from \cite{FUfindiff19} that are listed in {\tt Polyksrestcoeff3.m} in the {\tt j\_s} format.\\[0.5mm]  
For other matrix optimization processes it is important to reformulate the code development process above  and  to try and understand the interaction between suitable $\eta$ and $\tau$ values for discretized ZNN methods here in order   to be able to use ZNN well, see Remark 3 (a) and (b) below.\\[1mm]
An introduction to constrained optimization methods is available at  \cite{SD06}; see also \cite{LSX20}. Several optimization problems are studied in \cite{LSX20} such as  Lagrange optimization for unconstrained time--varying convex nonlinear optimizations called  U-TVCNO  and time--varying linear inequality systems called TVLIS. The latter will  be treated in subpart {\bf (VI)} just below.\\[2mm]
 {\bf Remark 3 : } (a) The computed results of any Lagrangian optimization algorithm should always be carefully scrutinized against what is well known in optimization theory, see  \cite{Fi76, Ro76, Gau77, BNO03, NW06}  for analytic requirements. Do not accept your ZNN computed results blindly. Our code development here applies to just one specific time–varying discretized matrix. 
   \\[0.3mm]
(b) An important concept  in ZNN's realm is   the product $ \tau \cdot \eta$ of the sampling gap $\tau$ and the exponential error decrease constant $\eta$ for any one specific problem and any  discretized  ZNN method which uses a specific suitable  finite difference scheme of one fixed type {\tt j\_s}. This product of constants, regularly denoted as $h = \tau \cdot \eta$ in the Zeroing Neural Network literature, seems to be nearly constant for the optimal choice of the parameters $\tau$ and $\eta$ over a wide range of sampling gaps $\tau$ if the chosen difference formula of type {\tt j\_s}  stays fixed.  Yet the optimal value of the near 'constant' $h$ varies widely from one look--ahead convergent finite difference formula to another. The reason for this behavior is  unknown and worthy of further studies; see the optimal $\eta$ for varying sampling gaps $\tau$ tables for time--varying eigenvalue computation in \cite{FUaccFOVLama20}. It is interesting to note that analytic continuation ODE methods also deal with optimal bounds for a product, namely that of the step--size $\tau$ and the local Lipschitz constant $L_i$ of the solution $y$, see \cite[p. 429]{EMU96}.\\[-0mm]
 \hspace*{35mm} \underline{\hspace*{90mm}}\\[-2mm]
 
 Thus far in this practical section, we have worked through five models and a variety of 
 time--varying matrix problems. We have developed seven detailed Matlab codes.  Our numerical codes all implement  ZNN look--ahead convergent difference formula based discretized processes  for time--varying matrix and vector problems in seven steps as outlined in Section 1. Each of the resulting ZNN computations requires a linear equations solve (or an inverse matrix times vector product) and a simple convergent recursion per time step. Some of the codes are very involved such as for example {\bf (V)} which relies on Matlab's symbolic toolbox and its differentiation functionality. Others were straightforward. All of our seven algorithms are look--ahead and rely only on earlier and current data to predict  future solutions. They do so with high accuracy and  run in fractions of a second over sampling gap $\tau$ time intervals  that are 10 to 100 times longer than their CPU run time. This adaptability makes discretized ZNN methods highly useful for real--time and on--chip implementations.\\[1mm]
 We  continue with further examples and restrict our explanations of discretized ZNN methods to the essentials from now on. We also generally refrain from coding further ZNN based programs now, except for one code in subpart  {\bf (VII start-up)} that benefits from a Kronecker  product formulation and where we also explain how to generate start--up data  from completely random first settings via ZNN steps.  We include extended references and encourage our readers to try and implement their own ZNN Matlab codes for their specific time--varying matrix or vector problems along the lines of our detailed examples {\bf (I)} ... {\bf (V)} above and the examples {\bf (VI)} ... {\bf (X)} that follow below.\\[-2mm]
 
{\bf (VI) \ Time--varying Linear Equations with Linear Equation and Inequality Constraints :}\\[2mm]
We consider two types of linear  equation and linear inequality constraints here :\\[-6mm]
\begin{eqnarray*}
\hspace*{-0mm} \circledzero\!\!\!_a \hspace*{2mm} {\bf (A)} & & \ \text{ and \ \ \ \ \ \ \ \ \circledzero\!\!\!}_c \ \  {\bf (AC)}\\
& \hspace*{-6mm} A(t)_{m,n} x(t)_n \ \leq \ b(t)_m & \hspace*{26mm}  \ A(t)_{m,n} x(t)_n \ = \  b(t)_m \hspace*{16mm}\\
&& \hspace*{26.2mm} \ C(t)_{k,n} x(t)_n \ \leq\ d(t)_k \ \ .
\end{eqnarray*}

\vspace*{-2mm}
We assume that the matrices and vectors all have real entries and that the given inequality problem has a unique solution for all $t_o \leq t \leq t_f \subset \RR$. Otherwise with subintervals of $[t_o,t_f]$ in which  a given problem is unsolvable or has infinitely many solutions, the problem itself would become subject to potential bifurcations and  thus well beyond the scope of this introductory ZNN survey paper.\\[1mm]
In general, solving parametric NLPs with inequality constraints is a difficult problem. The `active set'  -- the set of inequality constraints at their boundaries  --  can change, there can be so-called weakly-active constraints where a constraint is active but its associated multiplier is zero or the Lagrange multipliers may not be unique (albeit contained within a bounded convex polytope). For example, \cite{RD95} addresses the problem of non-unique multipliers by solving a linear program. To avoid such difficulties here, we introduce the idea of `squared slack variables', see, e.g.\ \cite{SW70},  to replace the  given  system with linear inequalities by a system of linear equations. This is usually a dangerous move due to the potential for numerical instabilities, e.g.\ \cite{Ro76, FF17}, and arguably should not be used in practice. It may suffice as a motivating example.\\[2mm]
 The slack variable vector $u^{\texttt{.}2} \in \RR^\ell$ typically has non-negative entries  in the form of real number squares, i.e., $u^{\texttt{.}2}(t) = {\tt [(u_1(t))^2; ~...~ ;  (u_\ell(t))^2]} \in \RR^\ell$ for $u(t) = {\tt [u_1;~ ... ~; u_\ell]}$ in Matlab column vector notation with $\ell = m$ or $k$, depending on the type of time--varying inequality system {\bf (A)} or {\bf (AC)}.  With $u$ our models {\bf (A)} and {\bf (AC)}  become  \\[-8mm]

\begin{eqnarray*}
\hspace*{-2mm}\circledzero\!\!\!_a \ \text{ for } {\bf (Au)} & &  \hspace*{2mm}\text{ and } \ \circledzero\!\!\!_c \ \text{ for } \  {\bf (ACu)}\\
 &\hspace*{-16mm}A(t)_{m,n} x(t)_n + u(t)_m = ~ b(t)_m  
 & \hspace*{4mm} \bp A_{m,n}& O_{m,k}\\  C_{k,n} & \text {diag}(u) \ep_{m+k,n+k} \bp x_n\\ u \ep_{n+k}  \ = \  \bp b_m\\ d_k \ep    \in \RR^{m+k} .
\end{eqnarray*}

\vspace*{-1mm}
For {\bf (ACu)},  diag$(u)$ denotes the $k$ by $k$ diagonal matrix with the entries of $u \in \RR^k$ on its diagonal.\\[.8mm]
 The error function for {\bf (Au)} and $u \in \RR^k$ is $\circledone\!\!\!_a  \ \ E(t) = A(t) x(t) + u(t) - b(t)$ . The product rule of differentiation, applied  to each component function of $u(t)$ establishes the error function DE for {\bf (Au)} as
 
 \newpage
 
 \vspace*{-5mm}
$$ \circledtwo\!\!\!_a\ \ \   \dot E = \dot A x + A \dot x + 2 ~ u~  .\!* \dot u - \dot b \ \stackrel{(*)}{=} \ -\et~ (A x + u - b) \ = \ -\et~ E  $$
where the $.*$ product uses the Matlab notation for entry--wise vector multiplication. If the unknown entries of $x\in \RR^n$ and $u\in \RR^m$ are gathered in one extended column vector $[x(t);u(t)] \in \RR^{n+m}$ we obtain the alternate error function DE for {\bf (Au)} in block matrix form as
$$ \circledtwo\!\!\!_a \ \ \ \dot E = \dot A x + \bp A&  2 \text{ diag}(u) \ep_{m,2m} \bp \dot x\\\dot u\ep_{2m} - \dot b \ \stackrel{(*)}{=} \ -\et~ (A x + v - b) \ = \ -\et~ E \ \in \RR^m \ .
$$
Similarly for  {\bf (ACu)},  the error function is $$\circledone\!\!_c \ \ \ E_c(t) = \bp  A & O \\  C & \text{ diag}(u) \ep_{m+k,n+k} \bp  x\\ u \ep_{n+k} - \bp  b\\  d \ep \in \RR^{m+k} $$
 and its error function DE is \circledtwo$\!\!_c$ \ 
$$\dot E_c(t) = \bp \dot A & O \\ \dot C & 2 \text{ diag}(u) \ep \bp \dot x\\\dot u \ep - \bp \dot b\\ \dot d \ep ~ \stackrel{(*)}{=}  ~ -\et~ \left( \bp A& O\\  C & \text {diag}(u) \ep \bp x\\ u \ep  -  \bp b\\ d \ep \right) =  -\et~ E_c(t) .
$$
Solving the error function DEs \circledtwo$\!\!_a$ and \circledtwo$\!\!_c$  for the derivative vectors  $[\dot x(t);\dot u(t)]$, respectively, via the built--in pseudo--inverse function {\tt pinv.m} of Matlab for example, see subsection ({\bf III}) above, we obtain the following expressions for the derivative of the unknown vectors $x(t)$ and $u(t)$.\\
For  model {\bf (Au)} 
$$ \circledthree\!\!_a \ \ \ \bp \dot x\\ \dot u \ep = ~ \text{pinv}\left((A_{m,n}\ \ \ 2~ \text{diag}(u) )_{m,n+m}\right)_{n+m,m} \cdot  \left( \dot b - \dot A x -\et~ (A x + u - b) \right)_m   \in \RR^{n+m} \ ,
$$
and for {\bf (ACu)} 
$$\circledthree\!\!_c \ \ \bp \dot x\\ \dot u\ep = ~ \text{pinv}\bp \dot A & O \\ \dot C & 2 \text{ diag}(u) \ep  \cdot  \left(
\bp \dot b\\ \dot d \ep  \ \  - \ \  \et~ \left( \bp A& O\\  C & \text {diag}(u) \ep \bp x\\ u\ep  -  \bp b\\ d \ep \right) \right) \in \RR^{n+k}  .
$$
The Matlab function {\tt pinv.m} in \circledthree uses the Singular Value Decomposition (SVD). The derivative of the vector $[x(t); u(t)] $ can alternately be expressed in terms of Matlab's least square function {\tt  lsqminnorm.m} in
$$ \circledthree\!\!_{als} \ \  \bp \dot x\\ \dot u\ep ~ = ~ \text{lsqminnorm}\left(\! ((A_{m,n}\ \ \ 2~ \text{diag}(u) )_{m,n+m})_{n+m,m} ,  \left( \dot b - \dot A x -\et~ (A x + u - b) \right)\! \right)_m   \in \RR^{n+m} \ ,
$$
or  
$$\circledthree\!\!_{cls}  \ \ \bp \dot x\\ \dot u \ep ~ = ~ \text{lsqminnorm} \left(\!
\bp \dot A & O \\ \dot C & 2 \text{diag}(u) \ep\!,  \left(\!
\bp \dot b\\ \dot d \ep\!    -   \et~ \left(\! \bp A& O\\  C & \text {diag}(u) \ep\! \bp x\\ u \ep  -  \bp b\\ d \ep\! \right)\! \right)\! \right)  .
$$
Next choose a look--ahead finite difference formula of type {\tt j\_s} for the discretized problem  and equate its derivative  $[\dot x(t_k);\dot u(t_k)]$ with the above value in \circledthree$\!\!_a$,  \circledthree$\!\!_{als}$ or \circledthree$\!\!_c$, \circledthree$\!\!_{cls}$ in order to eliminate the derivatives from now on. Then solve the resulting solution--derivative free equation for the 1--step ahead unknown $[x(t_{k+1}); u(t_{k+1})]$ at time $t_{k+1}$.\\[1mm]
 The Matlab coding  of a ZNN based discretized algorithm for time--varying linear systems  with equation or inequality constraints  can now begin after $j+s$ initial values have been set. \\[1mm]
Recent work on discretized ZNN methods for time--varying matrix inequalities is available in \cite[TVLIS]{LSX20}, \cite{ZYYHXH19} and \cite{XLNSG19}.

\newpage

{\bf (VII) \ Square Roots of Time--varying Matrix flows :}\\[2mm]
Square roots $X_{n,n} \in \CC^{n,n}$ exist for all nonsingular static matrices $A \in \CC^{n,n}$, generalizing the fact that  all complex numbers have square roots over $\CC$. Like square roots of numbers,  matrix square roots may be real or complex. Apparently  the number of square roots that a given nonsingular matrix  $A_{n,n}$ with $n \gg 2$ has is not known, except that  there are many. Every different set of starting values for our otherwise identical  nonsingular matrix flow $A(t)$ example of Figure 1 did result in a different matrix square root matrix flow $X(t)$ for $A(t)$ when covering $0 \leq t \leq 360 \ sec$.\\
 For singular matrices $A$ the existence of square roots depends on $A$'s Jordan block structure and its  nilpotent Jordan blocks $J(0)$ and some matching  dimension conditions thereof, see e.g. \cite[l. 15-18, p. 466; Thm. 4.5, p. 469; and Cor. 11.3, p. 506]{EU92} or the result in 
 \cite[Theorem 2]{CL74} formulated in terms of the ``ascent sequence" of $A$ for details.\\[1mm] 
Here we assume that our time--varying flow matrices $A(t)$ are nonsingular for all $ t_o \leq t \leq t_f \subset \RR$. Our model equation is \circledzero $A(t) = X(t)\cdot X(t)$ for the unknown time--varying square root $X(t)$ of $A(t)$. Then the error function becomes $\circledone  E(t) = A(t) - X(t)\cdot X(t)$ and the error DE under exponential decay stipulation is
$$ \circledtwo \dot E = \dot A - \dot X X - X \dot X  \stackrel{(*)}{=} - \et~ (A - X X) = - \et~ E \ $$  
where we have again omitted the time variable $t$ for simplicity. Rearranging the central equation $(*)$ in  \ \circledtwo \ with all unknown $\dot X$ terms on the left-hand side  gives us 
$$ \circledthree \ \ \dot X X + X \dot X = \dot A + \et~ (A - X X)\ . $$
Equation \circledone  is  model (10.4)  in \cite[ch. 10]{ZHYLLQ20} except for a minus sign.
In \circledthree  we have a similar situation as was encountered in Section 1 with the $n$ by $n$ matrix eigenvalue equation  for finding the complete time--varying matrix eigen--data. Here again,  the unknown  matrix derivative  $\dot X$  appears as  both a left and right factor in matrix products.  In Section 1 we switched our model and solved the time--varying matrix eigenvalue problem one eigenvector and eigenvalue pair at a time. If we use the Kronecker product for matrices and column vectorized matrix representations  in \circledthree  -- we could have done the same  in  Section 1 for the complete time--varying matrix eigen--data problem  -- then this model  can be solved directly.  And we can continue with the global model and  discretized ZNN when relying on notions from classical matrix theory, i.e., static matrix theory  helps to  construct a   discretized  ZNN algorithm for  time--varying matrix square roots.\\[1mm]
For two real or complex matrices $A_{m,n}$ and $B_{r,s}$ of any size,  the {\bf \emph{ Kronecker product}} is defined as the matrix
$$ A \otimes B  = \bp a_{1,1}B & a_{1,2}B & \dots & a_{1,n}B\\
a_{2.1}B &\ddots& & a_{2,n}B\\
\vdots && \ddots & \vdots\\
a_{m,1}B& \dots&\dots&a_{m,n}B \ep_{m\cdot r,n \cdot s}   \ .  $$ 
The command {\tt kron(A,B)} in Matlab creates  $A \otimes B$ for general pairs of matrices. Compatibly sized Kronecker products are added entry by entry just as matrices are.
For matrix equations a very useful property of the Kronecker product is the  rule  
\begin{equation*} \hspace*{46mm}
(B^T \otimes A )X(:) = C(:) \ \ \ \text{ where } \ \ C = AXB\  .
\hspace*{42mm} (21)
\end{equation*}
Here the symbols $X(:)$ and $C(:) \in \CC^{m n}$  denote the column vector storage mode in Matlab of any matrices $X_{m,n}$ or $C_{m,n}$. 

When we combine the Kronecker product  with Matlab's  column vector matrix notation  $M(:)$  we can rewrite the left--hand side of equation \circledthree $\dot X X + X \dot X = \dot A + \et~ A - \et~ X X $ as
$$ \hspace*{39mm}(X^T(t) \otimes I_n + I_n \otimes X(t))_{n^2,n^2}\cdot  \dot X(t)(:)_{n^2,1} \in \CC^{n^2} \ , \hspace*{38mm} (22)$$
while its right--hand side translates into
$$ \dot A(t)(:)_{n^2,1} + \et\ A(t)(:)_{n^2,1} - \et\  (X^T(t) \otimes I_n)_{n^2,n^2} \cdot X(t)(:)_{n^2,1} \in \CC^{n^2} \ . $$
 And miraculously the difficulty of $\dot X$ appearing on both sides as a factor in the 
 left--hand  side matrix products in \circledthree is gone.  We generally cannot tell whether the sum of Kronecker products in front of $\dot X(t)(:)$ in $(22)$ is nonsingular. But if we assume it is, then we can solve \circledthree for $\dot X(t)(:)$.
\begin{equation*}  \hspace*{4mm}\circledthree\!\!\!_K \ \  \dot X(t)(:) = (X^T(t) \otimes I_n + I_n \otimes X(t))^{-1} \cdot \left( \dot A(t)(:) + \et\ A(t)(:) - \et\ (X^T(t) \otimes I_n) \cdot X(t)(:)  \right) \ . \hspace*{4mm}(23)
\end{equation*}
Otherwise if  $X^T(t) \otimes I_n + I_n \otimes X(t)$ is singular, we simply replace the matrix inverse by the pseudo--inverse\\
 pinv$(X^T(t) \otimes I_n + I_n \otimes X(t))$ above and likewise in the next few lines. A recent paper \cite{ZYGLZ21} has dealt with a singular Kronecker product that occurs when trying to compute the Cholesky decomposition of  a positive definite matrix via ZNN. With 
$$P(t)  = (X^T(t) \otimes I_n + I_n \otimes X(t)) \in \CC_{n^2,n^2}$$ 
when assuming non-singularity and 
$$ q(t) = \dot A(t)(:) + \et\ A(t)(:) - \et\ (X^T(t) \otimes I_n) \cdot X(t)(:)  \in \CC^{n^2}$$ 
we  now have $\circledthree \!\!\!_K \  \dot X(t)(:)_{n^2,1} = P(t)\backslash q(t)$. This formulation is  reminiscent of formula (3$i$)  in Step 3 of Section 2, except that here the entities followed by $(:)$ represent column vector matrices instead of square matrices and these vectors now have  $n^2$ entries instead of $n$ in (3$i$). This might lead to execution time problems for real-time applications if the size of the original system is in the hundreds or beyond, while $n = 10$ or $20$ should pose no problems at all. How to mitigate such size problems, see \cite{N51510} e.g..\\[1mm]
To obtain the derivatives $\dot X(t_k)$ for each discrete time $t_k = t_o + (k-1) \tau$ for use in Step 5 of discretized ZNN, we need to solve the $n^2$ by $n^2$ linear system $P(t)\backslash q(t)$ and obtain the column vectorized matrix $\dot X(t_k)(:)_{n^2,1}$. Then we reshape $\dot X(t)(:)_{n^2,1}$ into square matrix form via Matlab's {\tt reshape.m} function. Equation \circledfive then equates the above matrix version $\dot X(t_k)_{n,n}$ of step  \circledthree with  the difference formula for $\dot X(t_k)_{n,n}$  from our chosen finite difference expression in step  \circledfour of ZNN and this helps us to predict $X(t_{k+1})_{n,n}$ in step  \circledsix.  This  has been done many times before in this paper, but without the enlarged Kronecker product matrices and vectors  and it should create no problems for our readers. A plot and analysis of the error function's decay as time $t$ progresses for 6 minutes was given in Section 2, preceding Figure 1,  for  the time--varying matrix square root problem.\\[1mm]
 For further analyses, a convergence proof and numerical tests of ZNN based time--varying matrix square root algorithms see \cite{ZHYLLQ20}. Computing time--varying matrix square roots is also the subject of \cite[Chs. 8, 10]{ZGbook15}.\\[-2mm]

{\bf (VIII) \ Applications  of 1--parameter Matrix Flow Results to Solve a Static Matrix Problem :}\\[1mm]
Concepts and notions of classical matrix theory  help us often with time--varying matrix problems. The concepts and results of the time--varying matrix realm can likewise help  with classic, previously unsolvable fixed matrix theory problems and applications. Here is one example.  \\[1mm]
 Numerically the Francis QR eigenvalue algorithm 'diagonalizes' every square matrix  $A$ over $\CC$ in a backward stable manner. It does so for diagonalizable matrices as well as for derogatory matrices, regardless of their Jordan structure or of repeated eigenvalues. QR finds  a backward stable 'diagonalizing' eigenvector matrix similarity for any $A$. For matrix problems such as least squares, the SVD, or the field of values problem that are unitarily invariant,  classic  matrix theory does not know of any way to unitarily  block diagonalize  fixed entry matrices $A \in \CC_{n,n}$. If such block decompositions could be found computationally, unitarily invariant matrix problems could be decomposed into  subproblems and thereby speed up the computations for decomposable matrices $A$.\\[0.5mm]    
  An  idea that was inspired by studies of time--varying matrix eigencurves in \cite{FUEigencross} can be adapted to find unitary  block decompositions of  static matrices $A$.  \cite{FUDecompMatrFlow23} deals with general 1--parameter matrix flows $A(t) \in \CC_{n,n}$. If $X$ diagonalizes one very specific  associated flow matrix ${\cal F}_A(t_1)$ via a unitary similarity $X^* ... X$ and $X^* {\cal F}_A(t_2)X$ is properly block diagonal for some $t_2 \neq t_1$  with ${\cal F}_A(t_1) \neq {\cal F}_A(t_2)$, then every ${\cal F}_A(t)$ is simultaneously block diagonalized by $X$ and consequently the flow $A(t)$ decomposes uniformly. \cite{FUDecompMatrFoV} then applies this specific matrix flow result to the previously intractable field of values problem for decomposing matrices $A$ when using path following methods. Here are the details.\\
  For any fixed entry matrix $A \in \CC_{n,n}$ the   hermitean and skew parts
  $${H=  (A+A^*)/2 = H^*} \ \ \text {  and } \ \  {K=  (A - A^*)/(2i) =  K^*}  \in \CC_{n,n} $$
of $A$   generate the  1--parameter  hermitean matrix flow
 $$ {\cal F}_A(t) = \cos(t)  H +   \sin(t)   K  = ({\cal F}_A(t))^*  \in \CC_{n,n}$$  for  all angles $0 \leq t \leq 2\pi$. This matrix flow goes back to Bendixson \cite{Be} and  Johnson \cite{J78} who studied effective ways to compute matrix field of values boundary curves. Whether there are different associated matrix flows that could enable  such matrix computations  and what they might look like is an open problem. 
 If   two matrices ${\cal F}_A(t_1)$  and ${\cal F}_A(t_2)$ in this specific flow are properly block diagonalized simultaneously by the same unitary matrix $X$ into the same block diagonal pattern for some $A(t_2) \neq A(t_1)$, then every matrix ${\cal F}_A(t)$ of the flow ${\cal F}_A$ is uniformly block diagonalized by  $X$ and subsequently so is $A = H + iK$ itself, see \cite{FUDecompMatrFoV} for details.\\[1mm]
The matrix field of values (FOV) problem \cite{J78} is invariant under unitary similarities. The field of values boundary curve of any matrix $A$ can be determined by finding the extreme real eigenvalues for each hermitean ${\cal F}_A(t)$ with $ 0 \leq t \leq 2\pi$ and then evaluating certain eigenvector $A$--inner products to construct the FOV boundary points in the complex plane. One way to approximate the FOV boundary curve is to compute full eigenanalyses of hermitean matrices ${\cal F}_A(t_k)$  for a large set of angles $0 \leq t_k \leq 2\pi$ reliably via Francis QR as QR is a global method and  has no problems with 
eigen--decompositions of normal matrices. Speedier ways use path following methods such as initial value ODE solvers or discretized ZNN methods. But path following eigencurve methods cannot ensure that they find the extreme eigenvalues of ${\cal F}_A(t_k)$ if the  eigencurves of  ${\cal F}_A(t)$ cross in the interval $[0,2\pi]$. Eigencurve crossings  can only occur for unitarily decomposable  matrices $A$, see \cite{NW29}. Finding eigencurve crossings for decomposing matrices $A$ takes up a large part of \cite{LM18} and still fails to adapt IVP ODE  path following methods for all possible types of decompositions for  static matrices $A$.\\
 The elementary method of \cite{FUDecompMatrFlow23}  helps to solve the field of values problem for decomposable matrices $A$ for the first time without having to compute all eigenvalues  of each  hermitean ${\cal F}_A(t_k)$ by -- for example -- using  the global Francis QR algorithm. Our combined matrix decomposition and discretized ZNN method is up to 4 times faster than the Francis QR based global field of values method or any other IVP ODE analytic continuation  method. It  depicts the FOV boundary accurately and quickly for all decomposing and indecomposible matrices $A\in \CC_{n,n}$, see \cite{FUDecompMatrFoV} for more details and ZNN Matlab codes.\\[2mm]
{\bf (IX) \ Time--varying Sylvester and Lyapunov Matrix Equations :}\\[2mm]
{\bf (S) } The static {\bf \emph{Sylvester equation}} model 
$$\hspace*{-4mm} \circledzero \ \ \   A X + X B = C$$
with $A \in \CC_{n,n}$, $B \in \CC_{m,m}$, $C \in \CC_{n,m}$ is solvable for $X \in \CC_{n,m}$ if $A$ and $B$ have no common eigenvalues.\\[1mm]
 From the error function \circledone 
$  E(t) = A(t) X(t)+ X(t) B(t) - C(t)$  we construct the exponential decay error DE\  $\dot E(t) = -\et~ E(t)$ for a positive decay constant $\et$ and obtain the equation
$$\!  \circledtwo \ \ \dot E(t) = \dot A(t) X(t)+ A(t) \dot X(t)+ \dot X(t) B(t) +  X(t) \dot B(t) - \dot C(t)  \stackrel{(*)}{=}-\et~ (A(t) X(t)+ X(t) B(t) - C(t))  = -\et~ E(t) $$ 
  and  upon reordering the terms in $(*)$ we have
$$\hspace*{-4mm} \circledthree \ \ \   A \dot X+ \dot X B   = 
-\et~ (A X+ X B - C) - \dot A X -  X \dot B   + \dot C   ,$$
where we have  dropped all references to the time parameter $t$ to simplify reading. Using the properties of Kronecker products and column vector matrix representations as introduced in Section 2 ({\bf VII}) above we rewrite the left--hand side of \circledthree as
$$  (I_m \otimes A(t) + B^T(t) \otimes I_n)_{n\cdot m,n\cdot m} \cdot \dot X(t)(:)_{n\cdot m,1}  =  M(t) \dot X(t)(:) \in \CC^{nm}$$ 
and the right--hand side as 
$$ q(t) =  - (I_m \otimes \dot A + \dot B^T\otimes I_n)_{nm,nm}\cdot X(:)_{nm,1} + \dot C(:)_{nm,1} - \et~ ( (I_m \otimes  A + B^T\otimes I_n)\cdot X(:) 
 - C(:))_{nm,1} \in \CC^{n m} \ .$$
The  Kronecker matrix product  is  necessary here to express the two sided appearances of $\dot X(t)$ on the left--hand side of \circledthree$\!\!$. The right--hand side of \circledthree can  be expressed more simply in column vector matrix notation without using Kronecker matrices as  
  $$q(t) =  - ((\dot A \cdot X)(:)+ (X\cdot \dot B)(:))_{nm,1}  +\dot C(:)_{nm,1} - \et~ \left( (A\cdot X)(:) + (X\cdot B)(:) - C(:)\right)_{nm,1} \in \CC^{n m} \ .$$ 
Expressions such as  $(A\cdot X)(:)$ above denote the column  vector representation of  the matrix product $A(t) \cdot X(t)$. Thus we obtain the   linear system $ M(t) \dot X(t) = q(t)$ for $\dot X(t)$  in \circledthree with $M(t) = (I_m \otimes A(t) + B^T(t) \otimes I_n)   \in \CC_{n m, n m}$ when using either form of $q(t)$. And $\dot X(t)(:) \in \CC^{n m}$ can be expressed in various forms, depending on the solution method and the case of (non)--singularity of $M(t)$ as  $\dot X(t)(:) = M(t)\backslash q(t)$, 
$\dot X(t)(:) = \text{inv}(M(t)) * q(t)$, or\\[-2mm] 
$$   \dot X(t)(:) = \text{inv}(M(t)) \cdot q(t) \ \   \text{ or }  \ \  \dot X(t)(:) = \text{pinv}(M(t)) \cdot q(t)  \ \ \text{  or } \ \ \dot X(t)(:) = \text{lsqminnorm}(M(t)) , q(t)) \ ,$$
with the latter two formulations to be used in case  $M(t)$ is singular. \\[0.5mm]
Which form of $q(t)$ gives faster or more accurate results for $\dot X(t)(:)$ can be tested in Matlab by opening the {\tt $>\!\!>$ profile viewer} and running the discretized ZNN method for both versions of $q(t)$ and the various versions  of $\dot  X(t)$. We have also mentioned several methods in Matlab  to (pseudo--)solve  linear systems with a singular system matrices  such as $M(t)$ above. Users can experiment and  learn how to optimize  such Matlab codes for their specific problems and for the specific version of Matlab that is being used.\\
   Once $\dot X(t)(:)$ has been found in column vector form it must be reshaped in Matlab into an $n$ by $m$ matrix $\dot X(t)$. Next we have to equate  our computed derivative matrix $\dot X_k$ in the discretized version at time $t_k$ with a specific look--ahead finite difference formula expression for $\dot X_k$ in step \circledfive$\!\!$. The resulting 
   solution--derivative free equation finally is solved for the future solution $X_{k+1}$ of the time--varying Sylvester equation in step \circledsix of our standard procedures list. Iteration then concludes the ZNN algorithm for Sylvester problems.\\[2mm]
{\bf (L) }  A suitable  time--varying {\bf \emph{Lyapunov equation}} model is $$  \circledzero \ \   A(t)X(t)A^*(t) - X(t) + Q(t) = O_{n,n} \ \ \text{ with a hermitean flow } \ Q(t) = Q^*(t) . $$
  Its error function is  
$$\circledone  \ \ \ E(t) = A(t)X(t)A^*(t) - X(t) + Q(t) \ \ \ \stackrel{!}{=} O_{n,n} \in \CC_{n,n} \ .$$
Here all matrices are complex and square of size $n$ by $n$. Now we introduce  a shortcut and convert the matrix error equation \circledone immediately to its  column vector matrix with Kronecker matrix product  form
$$ \circledone \!\!\!_{(cvK)} \ \ \ E_K(:) = (\bar A \otimes A) X(:) - X(:) + Q(:) \in \CC^{n^2}$$ 
where we have used the formula  
$(AXA^*)(:) = (\bar A \otimes  A)X(:)$ and dropped all mention of  dependencies on $t$ for simplicity. Working towards the exponentially decaying differential error equation for $E_K(:)$, we note that derivatives of time--varying Kronecker products $U(t)\otimes V(t)$ follow the  product rule of differentiation
$$ \dfrac{\partial (U(t)\otimes V(t))}{\partial t} = \dfrac {\partial U(t)}{dt} \otimes V(t) + U(t) \otimes \dfrac{ \partial V(t)}{dt} $$
according to \cite[p. 486 - 489]{MN85}. 
With this shortcut to column vector matrix  representation, the  derivative of  the error function  $\circledone\!\!\!_{(cvK)}$ for $E_K(:)$ is 
$$  \dot E_K(:) = ((\dot {\bar A} \otimes A) +(\bar A \otimes  \dot A))X(:) + (\bar A \otimes  A))\dot X(:)  - \dot X(:) + \dot Q(:) \in \CC^{n^2} $$
And the error function DE  \ $ \dot E_K(:)  = - \et~ E_K(:)$ becomes
\begin{eqnarray*}
\circledtwo \ \ \ \dot E_K(:) &=&  (\bar A \otimes \dot A -I_{n^2})\dot X(:) + (\dot {\bar A} \otimes A +\bar A \otimes  \dot A)X(:)  + \dot Q(:) \ \   \\
  &\stackrel{(*)}{=}& -\et~ (\bar A\otimes A) X(:) + \et~ X(:) -\et~ Q(:) \ \ = \ \ -\et~ E_K(:)  \ .
\end{eqnarray*}
Upon reordering the terms of the central equation $(*)$ in \circledtwo we have the following linear system for the unknown column vector matrix $\dot X(:)$
$$ \circledthree \ \ \ (I_{n^2} - \bar A \otimes \dot A) \dot X(:) = (\dot{\bar A} \otimes A + \bar A \otimes \dot A) X(:) - \et~ (I_{n^2} - \bar A \otimes A)X(:) + \et~ Q(:) + \dot Q(:) \in \CC^{n^2} 
$$
where $\bar A$ is the complex conjugate matrix of $A$. For $M(t) = (I_{n^2} - \bar A \otimes A)  \in \CC_{n^2,n^2}$ \ and 
$$q(t)(:) = (\dot{\overline {A(t)}} \otimes A(t) + \overline {A(t)} \otimes \dot A(t)) X(t)(:) - \et~ (I_{n^2} - \overline {A(t)} \otimes A(t))X(t)(:) + \et~ Q(t)(:) + \dot Q(t)(:) \in \CC^{n^2} $$
 we have to solve the system
$ M(t) \dot X(t)(:) = q(t)$ for $\dot X(t)(:) \in \CC^{n^2}$ as was  explained earlier for the Sylvester equation.  In Step 5 we equate the matrix--reshaped expressions for $\dot X$ in \circledthree and the  chosen look--ahead convergent finite difference scheme expression for $\dot X$ from Step 4. Then we solve the resulting solution--derivative free equation for $X_{k+1}$ in Step 6 for discrete data predictively and thereby obtain the discrete time ZNN iteration formula.  These steps, written out in Matlab commands, give us the computer code for Lyapunov.\\[1mm]
Introducing Kronecker products and column vector matrix notations early in the construction of discrete ZNN algorithms is a significant short--cut for solving matrix equations whose unknown solution  $X(t)$ will occur in several different positions of time--varying matrix products. This is a simple new technique that  speeds up  discretized ZNN algorithm developments for such time--varying matrix  equation problems.\\[1mm]
 ZNN methods for  time--arying Sylvester equations have recently been studied in \cite{XZDLL20}. For a new  right and left 2--factor  version of Sylvester see \cite{ZLLYH20}. For recent work on discretized ZNN and Lyapunov, see \cite{SL20} e.g.. \\[-1mm]

{\bf (X) \ Time--varying Matrices, ZNN Methods and Computer Science : }\\[2mm]
The recent development of new algorithms for time--varying matrix applications  has implications for our understanding of computer science and of tiered logical equivalences on several levels in our mathematical realms.\\
 The most stringent realm of math is `pure mathematics' where theorems are proved and where, for example, a square matrix either has a determinant equal to 0 or it has not.\\
  In the next, the mathematical computations realm with its floating point arithmetic, zero is generally computed inaccurately as not being 0 and  any computed value with magnitude below a threshold such as a small multiple or a fraction of the machine constant $eps$ may be treated rightfully as 0. In the computational realm the aim is to approximate quantities to high precision, including zero and never worrying about zero exactly being 0.\\
    A third, the least stringent realm of mathematics belongs to the engineering world. There one needs to find solutions that  are good enough to approach the ``true theoretical solution" of a problem as known from the `pure' realm asymptotically; needing possibly only 4 to 5 or 6 correct leading  digits for a successful algorithm.\\[1mm]
The concept of differing math--logical equivalences in these three tiers of mathematics is exemplified and interpreted in Yunong Zhang and his research group's recent paper \cite{ZYQZ20} that is well worth reading and contemplating about.\\[-3mm]
 \hspace*{35mm} \underline{\hspace*{90mm}}\\[-1mm]
 
{\bf (VII start-up) \ How to Implement ZNN Methods On--chip  for Sensor Data Inputs and Remote Robots :}\\[2mm]
 Let $A(t_k)$ denote the sensor output that arrives at time $t_o \leq t_k = t_o + (k-1) \tau \leq t_f$ for a given time--varying matrix `problem'. Assume further that this data arrives with the standard clock rate of 50 Hz at the on--board  chip  of a robot and that there is no speedy access to software such as Matlab as the robot itself may be autonomously running on Mars. At time instance $t_k$ the `problem' needs to be solved on--chip predictively in real-t-ime, well before time $t_{k+1}$. We must predict or compute  the problem's solution  $ x(t_{k+1})$ or $X(t_{k+1})$ on-chip with its limited resources and  solve the given problem (at least approximately) before   $t_{k+1}$ arrives.\\[2mm]
   How can  we generate start--up data for a discretized  ZNN method and the relatively large constant  sensor sampling gap  $\tau = 0.02\ \text{sec} = 1/50\ \text{sec}$ that is standardly used  in real--world applications. How can one   create start--up data  \emph{ out of thin air}. For real-time sensor data flows $A(t)$, we  assume that the robot has no information  about the  `theoretical solution' or that there may be no known `theoretical solution' at all. Without any usable a priori system information,   we have to construct the $j + s$ initial values for the unknown $x(t_k)$ or $X(t_k)$ and  $k \leq j+s$ that will then be used to iterate with a {\tt j\_s} look--ahead convergent finite difference scheme based discrete ZNN method.\\[1mm]
    After many tries at this {\em task}, our best choice for the first 'solution'  turned out to be a random entry  vector or matrix for $x(t_o)$ or $X(t_o)$ of proper dimensions and then iterating  through $j + s$ simple low truncation error order ZNN steps until a higher truncation order predictive ZNN method can take over for times $t_{k+1}$, running more accurately with local truncation error order $O(\tau^{j+2})$ when $k > j+s$.\\[1mm]
 Here we illustrate this random entries start--up process for the time--varying matrix square root example of solving $A(t)_{n,n} = X(t) \cdot X(t) \in \CC_{n,n}$ from subpart  {\bf (VII)} by using \\[1mm]
 (a) the Kronecker form representation of $\dot X(t) (:)$  of (23),  denoted here as (23$v$) for conformity,\\
 (b) the matrix column vector notation $X(:) \in \CC^{n^2}$ for $n$ by $n$ solution matrices $X$, and \\[0.5mm]
 (c) the Kronecker product rule  (21) \
  $(B^T \otimes A )X(:) = C(:) $ \
   for compatible matrix triple  products $ C = AXB$. \\[1mm]
   The Kronecker product representation requires two different notations of matrices $X$ here: one as a square array, denoted by the letter $m$ affixed to the matrix name such as in $X\!m$, and another as a column vector matrix, denoted by an added $v$ to the matrix name as in $X\!v$. With these notations, equation (23) now reads as
 $$\hspace*{28mm}  \dot X\!v = (X\!m^T \otimes I_n + I_n \otimes X\!m)^{-1} \cdot \left( \dot Av+ \et\ (Av -  (I_n \otimes X\!m) \cdot X\!v ) \right)  \hspace*{25mm} (23v)
 $$ 
 where we have again  dropped all mentions of the time parameter $t$ for ease of reading.\\[1mm]
 As convergent look-ahead finite difference formula at start--up we use the simple Euler rule of type {\tt j\_s} = {\tt 1\_\!\_2} 
 \begin{equation*}
\hspace*{60mm}\dot x(t_k) = \dfrac{x(t_{k+1} ) - x(t_k)}{\tau}
\hspace*{56mm} (24)
\end{equation*} 
which gives us the expression  $ x(t_{k+1}) = x(t_k) + \tau \dot x(t_k)$. When applied to the solution matrix $X(t)$ and combined with (23$v$), we obtain the explicit start--up iteration rule
$$
X\!v(t_{k+1}) = X\!v(t_k) + \tau \cdot  (X\!m^T \otimes I_n + I_n \otimes X\!m)^{-1} \cdot \left( \dot Av(t_k) + \et\ (Av(t_k)  -  (I_n \otimes X\!m(t_k)) \cdot X\!v(t_k) ) \right) \ 
$$
where $(I_n \otimes X\!m(t_k)) \cdot X\!v(t_k)$ expresses the matrix square $X\!m(t_k) \cdot X\!m(t_k)$ according to the Kronecker triple product rule (21) by observing  that $X\!m (t_k)\cdot X\!m(t_k) \cdot I = (I \otimes X\!m(t_k)) \cdot X\!v(t_k)$. Every iteration  step in the discretized ZNN algorithm is coded exactly like this equation as done  many times before, with time--adjusted expressions of $\dot Av$ and for a different look-ahead convergent finite difference formula of type {\tt j\_s} = {\tt 4\_5} accordingly in the main iterations phase when $k > j+s = 9$. \\[1mm]
The Matlab m-file {\tt  tvMatrSquareRootwEulerStartv.m} for finding time--varying matrix square roots `from scratch'  is available together with two auxiliary m-files in \cite{FUZNNSurvey}. \\[2mm]
To plot the error graph in Figure 1 of Section 2 with our code we compute 18,000 
time--varying  matrix square roots predictively for 360 {\em sec} or 6 minutes, or one every 50th of a second. This process takes around 1.3 seconds of CPU time which equates to 0.00007 seconds for each square root computation step and shows that  our discretized ZNN method is very feasible to run predictively in real--time for matrix square roots during each of the 0.02 second  sampling gap intervals. \\[1mm]
Note that any ZNN method that is based on a finite difference scheme {\tt  j\_s} such as  {\tt 4\_\!\_5} with $j = 4$ has a local truncation error order of $O(\tau^{j+2}) = O(0.02^{~6}) \approx O(6.4 \cdot 10^{-11})$. In this  example the model equation $A(t) = X(t) \cdot X(t)$ is satisfied from a random Euler based start--up after approximately 20 seconds with smaller than $10^{-10}$ relative errors. The model's relative  errors  decrease to around $10^{-13}$ after about 6 minutes according to Figure 1.\\[1mm]
As $\tau = 0.02$  {\em sec} is  the standard fixed 50 Hz clocking cycle for sensor based output, it is important to adjust the decay constant $\et$ appropriately for convergence: If a look--ahead ZNN method  diverges for your problem and one chosen $\eta$, reduce $\eta$. If the error curve decays at first, but has large variations after the early start--up phase and 10 or 20 seconds have passed, increase $\eta$. The optimal setting of $\eta$ depends on the data input $A(t_k)$ and its variations size, as well as the chosen finite difference formula. The optimal setting for $\eta$ cannot be predicted.\\
 This is one of the many open problems with time--varying matrix problems and both discretized and continuous ZNN methods that needs a deeper  understanding of the numerical analysis of time--varying matrix computations.\\[1mm]
Eight different model equations for the time--varying matrix square root problem are detailed in \cite[Ch. 10]{ZGbook15}. 
\\[-7mm]

\section{Conclusions}
This paper has  tried to explain the inner workings and describe the computational phenomena  of   a recent, possibly new or slightly different branch of numerical analysis for discretized  time--varying matrix systems from its inside out. Time--varying discretized ZNN matrix algorithms  are  built from standard concepts and well known relations and facts of algebra, matrix theory, and also aspects of static  numerical matrix analysis. Yet they differ in execution  from their near cousins,  the analytic continuation IVP differential  equation  solvers; in speed, accuracy, predictive behavior and more. Zhang Neural Networks do not follow the modern call for  backward stable computations that find  the exact solution of a nearby problem  whose distance from the given problem  depends on the problem's conditioning.  Instead Zhang Neural Networks compute highly accurate future solutions  based on an exponentially decaying error function that ensures their  global convergence from nearly arbitrary start-up data. In their  coded versions, the thirteen discretized ZNN time--varying matrix algorithm examples in this paper use just one linear equations solve  and a short recursion of  earlier systems data per time step, besides  some auxiliary set-up et cetera functions. These codes run extremely fast using previous data immediately  after time $t_k$ and they arrive  well before time $t_{k+1}$ at an accurate prediction of the unknown  variable(s) at the next time instance $t_{k+1}$.\\[1mm]
The standard version of discretized ZNN methods for time--varying matrix problems proceeds in seven steps. The ZNN steps do mimic standard IVP ODE continuation methods and it is challenging to try and understand the differences between them from the `outside'.  Here many time--varying problems from matrix theory and  optimization have been built from the ground up for ZNN, including working Matlab codes for most of them. This was done  with ever increasing levels of difficulty and complexity, from simple time--varying linear equations solving routines to time--varying matrix inversion and time--varying pseudo-inverses; from time--varying Lagrange multipliers for function optimization to more complicated matrix problems such as  time--varying matrix  eigenvalue problems, time--varying linear equations with inequality constraints, time--varying matrix square roots and  time--varying Sylvester and Lyapunov equations. Some of these algorithms require Kronecker product matrix representations and all can be built and handled  successfully in the discretized ZNN standard seven steps way. \\[1mm] 
On the way we have encountered  models and error functions for ZNN that do not yield usable derivative information for the unknown variable(s) and for which we have learned how to re-define our model and  its error function accordingly for success in ZNN.
We have pointed out alternate ways to solve the linear equations part of discretized ZNN differently in Matlab and shown how to optimize the speed of some ZNN Matlab codes. We have dealt with  simple random entry early on-chip  'solutions' for unknown and unpredictable sensor  data  as starting values of the discretized ZNN iterations for time--varying matrix square roots in light of Kronecker products.\\[1mm]
 But we have not dealt with  the final, the engineering application of discretized  time--varying sensor driven ZNN matrix algorithms or have we created on-chip circuit designs for use in the control of chemical plants and in robots, for autonomous vehicles,  and other machinery. That was the task of \cite{ZGbook15} and this is amply explained in Simulink schematics there. Besides, circuit diagrams  appear often in the  Chinese engineering literature that we have quoted.\\[1mm]
There are many wide open questions with ZNN as, for example, how to choose among dozens and dozens of otherwise equivalent same {\tt j\_s} type look-ahead and convergent finite difference formulas for improved accuracy or speed and also regarding the  ability  of specific finite difference formulas to handle widely varying sampling gaps $\tau$ well. Neither do we  know how to assess or distinguish between high optimal $h = \tau \cdot \et$ and low optimal $h$ value  finite difference formulas, nor why there are such variations in $h$ for equivalent truncation error order formulas. These are open  challenges for experts in difference equations.\\[1mm]
 Many other open questions with discretized ZNN are mentioned here and in some of the quoted ZNN papers.\\
 A rather simple test problem  would be to try and  solve the mass matrix differential equations \circledthree for the unknown solution $x(t)$ in Step 3 of some of our ZNN code developments with initial value ODE solvers such as {\tt ode23t}, {\tt ode45},  {\tt ode15s}, {\tt ode113} or {\tt ode23s}  in Matlab and see how well and quickly  the solution $x$ can be evaluated over time when compared with complete ZNN methods for the same problem. Such a comparison has appeared in \cite{FUaccFOVLama20} where our specific ZNN algorithm was compared to Loisel and Maxwell's \cite{LM18} directly differentiated eigendata equation and its IVP ODE solution.\\[-7mm]

\section{The Genesis of this Paper}
The author's  involvement with Zhang Neural Networks (ZNN) began in 2016  when he was sent  the book \cite {ZGbook15}, written  by Yunong Zhang and Dongsheng Guo from the Zentralblatt for review, see Zentralblatt 1339.65002. \\[1mm]
 The author of this introductory  survey  was impressed by the ideas and workings of ZNN. He  contacted Yunong Zhang and visited him and his research group at Sun Yat-Sen University in Guangzhou in the summer of 2017. Thanks to this visit and through  subsequent  exchanges of e-mails, questions  and advice back and forth, he began to enter  this totally new-to-him area of predictive time--varying matrix numerics. When writing and trying to submit papers on ZNN methods to western matrix  and  numerical and applied mathematics journals, he soon learned that  time--varying numerical matrix methods and predictive solutions thereof had been nearly untouched in the West; there were no suitable referees.  This  20 years old, yet new area had come to us from the East.  It resides and flourishes outside of our western  knowledge and understandings base,  with more than 90 $\%$ of its estimated 400  research and engineering papers and all of its at least five  books originating in China.   Only a few European and even fewer American scientists have  begun to contribute to the field, often  with publications that include Chinese coauthors; see, e.g., \cite{BA22, SWM19,UZTVeigenLAA19}. Moreover, in the emerging global ZNN engineering iterature of today there hardly is a hint of theoretical  work on ZNN and  no published numerical analysis for this new method and area.     \\[1mm]
In the summer of 2019 the author visited  Nick Trefethen and Yuji Nakatsukasa in Oxford. And we  conferred  for several hours about the numerical ideas that may lie  behind discretized Zhang Neural Networks and its error decay. \\[1mm]
In October 2020, when the first version of the manuscript was submitted, Nick Trefethen advised the editor to look for a referee with expertise in the vast literature of numerical methods for ODEs. And indeed, the group of  referees  discovered  analogies and possible connections between discretized ZNNs and certain analytic continuation methods, as described in Stephen Robinson's paper \cite{Ro76} from 1976, and in the books by Peter Deuflhard \cite{DNM11}, by  Eugene Allgower and Kurt George \cite{AG90}, and by Uri M. Ascher, Hongsheng Chin  and Sebastian Reich  \cite{ACR94},  and in the paper by  Uri M. Ascher and Linda R. Petzold \cite{AP98} that all originated in the 1990s.  More recent are the optimization books by Jorge Nocedal and Stephen J. Wright \cite{NW06} from 1999/2006  and by Dimitri Panteli Bertsekas, Angelia Nedi\'c and Asuman E. Ozdaglar \cite{BNO03} from 2003/2006,  and the paper by  Ellen Fukuda and Masao Fukushima \cite{FF17} from 2017. In  \cite{ACR94} the authors analyzed how to stabilize the DAEs in analytic continuations algorithms and introduced certain analytic continuation limit manifolds that might shed more understanding and light onto their connections with Zhang Neural Networks. But none of these expansive analytic continuation notions  in  \cite{DNM11} or \cite{AG90,ACR94,AP98} were known to Zhang and Wang in 2001 and no numerical analysis studies of their possible connections to ZNN have ever been attempted so far. The second round referee Peter Maxwell proved quite  knowledgable on these possibly adjacent method classes. Subsequently Maxwell and the author exchanged a series of extended reviews and rebuttals that eventually lead to the present form of the publication. Unfortunately we  could not settle the question of the level or actual kind of interconnectedness between  analytic continuation algorithms and the Zhang Neural Network method.\\[1mm] 
Eventually the seven step set-up structure of discretized Zhang Neural Networks became clear.  ZNN starts with a  global entry-wise error function for  time--varying matrix problems and  it stipulates an exponential error decay  that makes convergence and noise suppression automatic. Then ZNN continues  with an ingenious way to replace  the error function differential equation with  solving linear equations repeatedly instead and using high error order look-ahead convergent recursion formulas in tandem. Such difference formulas  had only been sparingly used in predictor-corrector ODE schemes before with relatively low error orders. To our knowledge, the discretized ZNN solution-derivative free setting has been achieved completely outside of the ODE analytic continuation realm and ZNN converges quite differently, see \cite[Figures 3 and 5]{FUAZNN23} for examples with convergence into the machine constant error levels (and below) over time. \\[1mm] 
This discovery and others such as the disparate term magnitudes of an adapted AZNN algorithm in  \cite[Table 1]{FUAZNN23} in the predictive computational Zhang Neural Network step have further strengthened this introductory and partial survey paper of Zhang Neural Network methods. As author I am glad and very thankful for the help and patience of my family and for the editor's and referees' helpful comments, to Nick and Yuri, to Peter  Benner who pointed me to  eigencurves  and  analytic continuation ODE methods (see \cite{FUaccFOVLama20}), and to Peter Maxwell  in particular for helping to improve the Lagrange optimization sections {\bf (V)} and  {\bf (VI)} in Section 3.\\[1mm]
 I do hope that time--varying matrix methods and continuous or discretized ZNN time--varying matrix methods can  become a new part of our global numerical matrix analysis research and that they will add to our  numerical matrix analysis know--how soon.\\[1mm]
 For this we  need to build a new time--varying matrix numerics knowledge base in the West; just as has been done so many times before for our own ever evolving  matrix computational needs, see \cite {FUEpochsMath18} for example.\\[-6mm]


\vspace*{4mm}

        \hspace*{20mm} at \ \ \ {\url {http://arxiv.org/abs/2008.02724}} \quad  August 10, 2022\\[2mm]
        
 \noindent
\centerline{ Tex file at \ \  ... .. /Box/local/latex/ZNNSurvey21/ZeroingNeuralNetworks5abc.tex \  \quad \today }\\       

\vspace*{2mm}

\noindent
 
2 image files :\\[3mm]
TVMatrSquareRoot5c.png\\
CuneiformYBC4652.png

\end{document}